\documentclass[twocolumn]{autart}    

\usepackage{graphicx}          

\usepackage{amsmath,amsfonts,amssymb,color} 
\usepackage{hyperref}
\hypersetup{
    colorlinks=true,
    linkcolor=blue,
    citecolor=red,
    filecolor=magenta,      
    urlcolor=cyan,
           }                                                                 
\usepackage{epsfig}
\usepackage{epsfig}
\usepackage{psfrag}
\usepackage{algorithm}
\usepackage{algpseudocode}
\usepackage{array}
\usepackage{epstopdf}
\usepackage{cite}
\usepackage{mathtools}
\usepackage[utf8]{inputenc}
\usepackage{subfig}
\usepackage{caption}
\usepackage{bm}

\DeclarePairedDelimiter\norm{\lVert}{\rVert}
\DeclareMathOperator*{\argmax}{arg\,max}

\newcommand{\bw}{\mathbf{w}}
\newcommand{\bx}{\mathbf{x}}
\newcommand{\by}{\mathbf{y}}
\newcommand{\bz}{\mathbf{z}}
\newcommand{\bu}{\mathbf{u}}
\newcommand{\bg}{\mathbf{g}}
\newcommand{\bv}{\mathbf{v}}
\newcommand{\X}{\mathcal{X}}

\newcommand{\BB}{\mathbb{B}}

\newcommand{\BS}{\mathbb{S}}
\newcommand{\M}{\mathcal{M}}

\newcommand{\A}{\mathcal{A}}

\newcommand{\CO}{\mathcal{O}}
\newcommand{\R}{\mathbb{R}}

\newcommand{\Z}{\mathbb{Z}}
\newcommand{\E}{\mathbb{E}}

\newcommand{\V}{\mathcal{V}}

\newcommand{\I}{\mathcal{I}}
\newcommand{\BI}{\mathbb{I}}

\begin{document}

\begin{frontmatter}

\title{Online Mixed Discrete and Continuous Optimization: Algorithms, Regret Analysis and Applications} 


\author[Hust]{Lintao Ye}\ead{yelintao93@hust.edu.cn},    
\author[Hust]{Ming Chi}\ead{chiming@hust.edu.cn},
\author[Hust]{Zhi-Wei Liu}\ead{zwliu@hust.edu.cn},
\author[NJUPT]{Xiaoling Wang}\ead{xiaolingwang@njupt.edu.cn},
\author[Pu]{Vijay Gupta}\ead{gupta869@purdue.edu},               

\address[Hust]{School of Artificial Intelligence and Automation, Huazhong University of Science and Technology, Wuhan 430074, China}  
\address[NJUPT]{College of Automation and College of Artificial Intelligence, Nanjing University of Posts and Telecommunications, Nanjing 210023, China} 
\address[Pu]{The Elmore Family School of Electrical and Computer Engineering, Purdue University, West Lafayette, IN 47907, USA}             


\begin{abstract}                          
We study an online mixed discrete and continuous optimization problem where a decision maker interacts with an unknown environment for a number of $T$ rounds. At each round, the decision maker needs to first jointly choose a discrete and a continuous actions and then receives a reward associated with the chosen actions. The goal for the decision maker is to maximize the accumulative reward after $T$ rounds. We propose algorithms to solve the online mixed discrete and continuous optimization problem and prove that the algorithms yield sublinear regret in $T$. We show that a wide range of applications in practice fit into the framework of the online mixed discrete and continuous optimization problem, and apply the proposed algorithms to solve these applications with regret guarantees. We validate our theoretical results with numerical experiments.
\end{abstract}

\end{frontmatter}

\section{Introduction}\label{sec:intro}
In online optimization problems, a decision maker needs to make a decision and receives a reward at each discrete point (i.e., round). The challenge is that the decision maker is interacting with an {\it unknown} environment, i.e., the decision must be made in an online manner before receiving the reward. The goal is to maximize the accumulative reward after a certain number of rounds. Online optimization problems have been widely studied in the literature and various algorithms have been proposed to solve the problem with theoretical performance guarantees (e.g., \cite{littlestone1994weighted,auer2002nonstochastic,zinkevich2003online,flaxman2005online,agarwal2010optimal,yang2016tracking}). 
However, most existing works focus on cases when the decision maker chooses a single discrete variable or a single continuous variable. For instance, in the Multi-Armed Bandit (MAB) problem, a decision maker needs to choose to play an arm (or a subset of arms) from all the candidate arms and receives a reward an associated reward (e.g., \cite{auer2002nonstochastic,lattimore2020bandit}). In the online convex (resp., concave) optimization problem, a decision maker needs to choose a continuous action (i.e., a point from a metric space) and receives the reward of the chosen action that is given by the value of a convex (resp., concave) function at that point (e.g., \cite{hazan2016introduction}).

In this paper, we propose to study an Online Mixed Discrete and Continuous Optimization (OMDCO) problem over a number of $T$ rounds, where the decision maker needs to choose a discrete variable and a continuous variable simultaneously at each round and receives a reward that depends on both the chosen discrete and continuous variables. The OMDCO problem generalizes the online optimization problems with a single discrete or continuous variable and captures a wider range of applications in practice. We summarize some related work as follows.

Online optimization with a single discrete variable has been widely studied, e.g., the MAB problem \cite{littlestone1994weighted,auer2002nonstochastic,bubeck2011introduction,lattimore2020bandit}, which is a benchmark setting for reinforcement learning \cite{sutton2018reinforcement}. Algorithms have been proposed to solve the MAB problem with sublinear regret in the total number of rounds in the problem (i.e., $T$), e.g., the {\bf Exp3} algorithm  \cite{auer2002nonstochastic}, which compares the accumulative reward of the algorithm after $T$ rounds and the optimal accumulative reward of a clairvoyant that knows the rewards of the arms in advance (e.g., \cite{auer2002nonstochastic}). In the classic MAB problem, the discrete variable is in fact a single element from a finite ground set. The formulation with the discrete variable being a subset chosen from the ground set  was studied in \cite{streeter2008online,golovin2014online,chen2016combinatorial}. Specifically, \cite{streeter2008online,golovin2014online} consider choosing the subset under a matroid constraint and propose an algorithm consisting of multiple parallel {\bf Exp3} subroutines, each of which chooses a single element to form the subset for each round. The algorithm achieves a sublinear approximate regret that compares the accumulative reward incurred by the algorithm and an approximately optimal accumulative reward. Differently, \cite{chen2016combinatorial} studies the combinatorial MAB framework with super arms, each of which is a subset of base arms from the ground set, and proposes an algorithm that selects a super arm in each round that achieves sublinear approximate regret.

Online optimization with a single continuous variable is also well-studied. In particular, when the reward of the continuous variable is given by a convex (or concave) function, various algorithms with provable sublinear regret in $T$ have been proposed to solve the problem (e.g., \cite{hazan2016introduction}). The proposed algorithms rely on extending the offline gradient-based algorithms to the online setting (e.g., \cite{flaxman2005online,agarwal2010optimal}).

Optimization problems with both discrete and continuous variables have also been studied, but in the offline setting where the objective function is fully known before the variables are chosen \cite{elenberg2018restricted,adibi2022minimax,bunton2022joint}. A well-studied instance is the mixed-integer program, where the objective function is typically assumed to be linear with a special structure \cite{richards2005mixed}. For general objective function, the authors in \cite{adibi2022minimax} proposed an iterative algorithm and showed that the algorithm converges to an approximately optimal solution.

We summarize the main challenges and contributions in this paper. First, we formulate the OMDCO problem. Since an algorithm for OMDCO needs to choose both discrete and continuous variables, the existing algorithms for MAB \cite{auer2002nonstochastic} and online concave optimization \cite{flaxman2005online,bartlett2007adaptive,agarwal2010optimal} that only deal with a single (discrete or continuous) variable cannot work. Thus, we design a novel algorithm and prove regret sublinear in $T$ under different settings. Our regret guarantees hold for the clairvoyant that can potentially switch its decision over the rounds in the OMDCO problem, which leads to the notion of dynamic regret that has been widely used to characterize the performance of algorithms for online optimization problems with a single discrete \cite{auer2002nonstochastic,garivier2011upper} or continuous variable \cite{jadbabaie2015online,li2020online}. The major challenge in analyzing the regret is that the rewards of discrete and continuous variables are coupled with each other. To tackle this, we carefully decompose the regret into two terms corresponding to the discrete and continuous variables, respectively. To bound the regret incurred by the discrete variable, we rely on extending the classical algorithms 
for the MAB problem to the setting when the decision maker can only receive an erroneous version of the reward of the arm played in each round. We also show that a wide range of important applications in practice fit into the general framework of OMDCO. To apply the algorithms proposed for the general OMDCO problem to these applications, we prove additional results and validate our theoretical results with numerical examples.

\paragraph*{Notation and Terminolog}
The sets of integers and real numbers are denoted as $\mathbb{Z}$ and $\mathbb{R}$, respectively. The set of integers (resp., real numbers) that are greater than or equal to $a\in\mathbb{R}$ is denoted as $\mathbb{Z}_{\ge a}$ (resp., $\R_{\ge a}$). The space of $m$-dimensional real vectors is denoted by $\mathbb{R}^{m}$ and the space of $m\times n$ real matrices is denoted by $\mathbb{R}^{m\times n}$. 
For a vector $\bx\in\R^n$, let $\bx^{\top}$ be its transpose and let $\norm{\bx}=\sqrt{\bx^{\top}\bx}$. For any $n\in\Z_{\ge1}$, let $[n]=\{1,\dots,n\}$. For any $\bx\in\R^n$ and any $S\subseteq[n]$, let $\bx(S)\in\R^n$ satisfy that $\text{supp}(\bx)=S$ and $(\bx(S))_i=(\bx)_i$ for all $i\in[n]$, where $(\bx)_i$ denotes the $i$-th element of the vector $\bx$ and $\text{supp}(\bx)\triangleq\{i\in[n]:(\bx)_i\neq0\}$. Let $\mathbf{1}$ be an all-one column vector whose dimension can be inferred from the context. For a real number $x\in\R$, let $|x|$ be its absolute value. For a finite set $S$, let $|S|$ be its cardinality. Let $\BB$ and $\BS$ denote the unit ball and unit sphere in $\R^n$ centered at the origin, respectively. For $\X\subseteq\R^d$ and $\xi\in\R$, let $\xi\X=\{\xi\bx:\bx\in\X\}$ and let $\Pi_{\X}(\bx)$ be the projection of $\bx\in\R^d$ onto $\X$.

\section{Problem Formulation and Preliminaries}\label{sec:problem formulation}
We study the problem of Online Mixed Discrete and Continuous Optimization (OMDCO) over $T\in\Z_{\ge1}$ rounds. At each round $t\in[T]$, the decision maker first chooses a discrete variable $S_t\in\I$ and a continuous variable $\bx_t\in\X\subseteq\R^d$, where $\I$ is a family of subsets of a finite set $[n]$. The decision maker then receives a reward given by $f_t(S,\bx)$, where $f_t:2^{[n]}\times\R^d\to\R_{\ge0}$. Let $(S_t,\bx_t)$ denote the decision point for any $t\in[T]$. To be more precise, we consider the scenario where the reward functions $f_1(\cdot,\cdot),\cdots,f_T(\cdot,\cdot)$ are not known to the decision maker a priori. One may view that the sequence of functions $f_1(\cdot,\cdot),\dots,f_T(\cdot,\cdot)$ are chosen by an adversary, and our algorithm design and regret analysis work for a potentially adaptive adversary \cite{bartlett2007adaptive,agarwal2010optimal}, i.e., the adversary is allowed to choose the function $f_t(\cdot,\cdot)$ based on the decision history $(S_1,\bx_1),\dots,(S_{t-1},\bx_{t-1})$. In this paper, we focus on the case when $\I$ is a uniform matroid, i.e., $\I=\{S\subseteq[n]:|S|\le H\}$ with $H\in\Z_{\ge1}$. In words, the decision maker is allowed to choose a subset from $[n]$ with cardinality at most $H$. We leave investigating the case when $\I$ is a general matroid \cite{lawler2001combinatorial} to future work. 

Our goal is to propose online algorithms for the decision maker and the performance of the online algorithms is characterized by the following $\alpha$-regret:
\begin{equation}\label{eqn:def of R alpha}
R(\alpha)\triangleq \E\Big[\alpha\big(\sum_{t=1}^T\max_{S\in\I,\bx\in\X}f_t(S,\bx)\big)-\sum_{t=1}^Tf_t(S_t,\bx_t)\Big],
\end{equation}
where the expectation is taken with respect to the randomness of the online algorithm, and $\alpha\in(0,1]$  is a parameter that will be specified later. In words, $R(\alpha)$ defined above measures the gap between the reward that the decision maker receives and an (approximately) optimal reward. Let us denote
\begin{equation}\label{eqn:optimal point}
(S^{\star}_t,\bx^{\star}_t)\in\argmax_{S\in\I,\bx\in\X}f_t(S,\bx),\forall t\in[T].
\end{equation}
The notion of regret has been widely used to characterize the performance of online algorithms with either a single discrete variable \cite{auer2002nonstochastic,streeter2008online,zhang2019online} or a single continuous variable \cite{flaxman2005online,bartlett2007adaptive,hazan2014bandit}. In particular,  we focus on the dynamic regret of online algorithms (e.g., \cite{auer2002nonstochastic,yang2016tracking,jadbabaie2015online}), i.e., the benchmark solution $(S_t^{\star},\bx^{\star}_t)$ defined in \eqref{eqn:optimal point} used in Eq.~\eqref{eqn:def of R alpha} can potentially vary over $t\in[T]$. A weaker notion of regret is termed as static regret where the benchmark solution cannot vary over $t\in[T]$ (e.g., \cite{zinkevich2003online}). For example, considering the benchmark solution $(S_1^{\star},\bx_1^{\star})=\cdots=(S_T^{\star},\bx_T^{\star})\in\argmax_{S\in\I,\bx\in\X}\sum_{t=1}^Tf_t(S,\bx)$ in Eq.~\eqref{eqn:def of R alpha}, the corresponding $R(\alpha)$ is called the static regret of the online algorithm. Since our analysis in this paper works for the stronger notion of regret (i.e., dynamic regret), it also holds if we consider static regret. To proceed, we introduce the following assumptions. 
\begin{assum}\label{ass:domains of f_t}
The set $\X\subseteq\R^d$ is convex, and there exist $r,D\in\R_{>0}$ such that $r\BB\subseteq \X\subseteq D\BB$.
\end{assum}
\begin{assum}\label{ass:function f_t}
For any $S\in\I$ with $|S|=H$ and any $t\in[T]$, $f_t(S,\cdot)$ is concave and $G$-Lipschitz over $\X$, and $f_t(S,\bx)\le C$ for all $\bx\in\X$, where $G,C\in\R_{\ge0}$.\footnote{A function $f:\R^d\to\R$ is $G$-Lipschitz over $\X\subseteq\R^d$ if $|f(\by_1)-f(\by_2)|\le G\norm{\by_1-\by_2}$ for all $\by_1,\by_2\in\X$.}
\end{assum}
\begin{assum}\label{ass:function f_t zero}
For any $t\in[T]$, $f_t(\cdot,\bx_t^{\star})$ is monotone nondecreasing with $f_t(\emptyset,\bx_t^{\star})=0$.\footnote{A set function $g:2^{[n]}\to\R$ is monotone nondecreasing if $g(S_1)\le g(S_2)$ for all $S_1\subseteq S_2\subseteq[n]$.} 
\end{assum}
Assumptions~\ref{ass:domains of f_t}-\ref{ass:function f_t} are made for the continuous part of the objective function, i.e., $f_t(S,\cdot)$, which are standard in online convex optimization (with a single continuous variable) \cite{zinkevich2003online,agarwal2010optimal,yang2016tracking}. Note that Assumption~\ref{ass:domains of f_t} requires $\mathbf{0}\in\X$, which is without loss of generality. 
Moreover, assuming $f_t(S,\bx)\in[0,C]$ for all $S\in\I$ with $|S|=H$ and $\bx\in\X$ is merely to ease our presentation. If $f_t(S,\bx)\in[f_a,f_b]$ for all $S\in\I$ with $|S|=H$ and all $\bx\in\X$, one can replace $0$ and $C$ with $f_a$ and $f_b$, respectively, in our analysis in the sequel.  Assumption~\ref{ass:function f_t zero} is made for the discrete part of the objective function, i.e., $f_t(\cdot,\bx_t^{\star})$, which is standard in combinatorial (or set function) optimization \cite{bian2017guarantees,streeter2008online,golovin2014online,das2018approximate}. Assumption~\ref{ass:function f_t zero} is a natural assumption, since adding more elements from the ground set $[n]$ to a set $S\in\I$ typically does not decrease the reward $f_t(S,\bx_t^{\star})$ of the decision point $(S,\bx^{\star}_t)$. We will justify this point later in   Section~\ref{sec:applications} using examples. Supposing the assumption $f_t(\emptyset,\bx_t^{\star})=0$ is removed, our analysis will work for the normalized function $f_t(S,\bx)-f_t(S,\bx^{\star}_t)$.

{\bf Partial information pattern.} We focus on two classes of information patterns: Single-point feedback and multi-point feedback. In single-point feedback (also known as the bandit information setting \cite{flaxman2005online,hazan2014bandit}), in any round $t\in[T]$, the decision maker observes only the function value of $f_t(\cdot,\cdot)$ at $(S_t,\bx_t)$ after choosing the point $(S_t,\bx_t)$ in round $t$. In multi-point feedback, in any round $t\in[T]$, the decision maker can observe the function value of $f_t(\cdot,\cdot)$ at multiple points, including $(S_t,\bx_t)$, after choosing the point $(S_t,\bx_t)$ in round $t\in[T]$. The single-point feedback is a natural setting in online optimization since only the value $f_t(S_t,\bx_t)$ of the decision point $(S_t,\bx_t)$ chosen in round $t\in[T]$ is revealed to the decision maker. Nonetheless, the evaluation of $f_t(\cdot,\cdot)$ at points other than the point $(S_t,\bx_t)$ chosen in round $t\in[T]$ may also be possible via prediction \cite{cesa2006prediction} or an external simulator \cite{xin2022identifying}. The single-point feedback and multi-point feedback are also known as partial information setting in online optimization. In contrast, in the full information setting \cite{zinkevich2003online}, the decision maker can directly observe the function $f_t(\cdot,\cdot)$. 
In this paper, we focus on the partial information case which is more challenging yet more practical and has been studied in online optimization with a single continuous variable \cite{agarwal2010optimal,yang2016tracking}.

\section{Results on Multi-Armed Bandit and Online Concave Optimization}\label{sec:preliminary results}
Our algorithm design and regret analysis rely on (variants of) the {\bf Exp3.S} algorithm for the MAB problem \cite{auer2002nonstochastic} and Online Concave Optimization (OCO) algorithms \cite{flaxman2005online,agarwal2010optimal}. We first provide some preliminary results on the MAB problem and OCO problem, which extend the results from \cite{auer2002nonstochastic,agarwal2010optimal,hazan2007logarithmic} and may be of independent interest.
\subsection{Multi-armed bandit with error feedback}\label{sec:multi-armed bandit}
\begin{algorithm}[tb]
\textbf{Input:} Parameter $\gamma\in(0,1]$
\caption{{\bf Exp3.S} with error feedback}\label{alg:exp3.S with error}
\begin{algorithmic}[1]
\State Set $\varpi_1^i=1$ for all $i\in[n]$
\For{$t=1$ to $T$}
\State Set $p^i_t=(1-\gamma)\frac{\varpi_t^i}{\sum_{j=1}^n\varpi_t^j}+\frac{\gamma}{n}$ for all $i\in[n]$
\State Choose $i_t\in[n]$ randomly accordingly to the probabilities $p_t^1,\dots,p_t^n$
\State Observe $\beta_t(r_t^{i_t}-\epsilon_t^{i_t}+\epsilon_t)$ 
    \For{$j=1$ to $n$}
    \State Set $\hat{r}^j_t\triangleq\begin{cases}\frac{\beta_t(r^j_t-\epsilon_t^j+\epsilon_t-a)}{p^j_t(b-a)}\ &\text{if}\ j=i_t\\0\ &\text{otherwise}\end{cases}$
    \State Set $\varpi_{t+1}^j=\varpi_{t}^j\text{exp}(\frac{\gamma\hat{r}_t^j}{n})+\frac{e}{nT}\sum_{j=1}^n\varpi_t^j$
    \EndFor
\EndFor
\end{algorithmic}
\end{algorithm}
An instance of the multi-armed bandit problem is given by a finite set $\I=[n]$ of possible actions (i.e., arms), and rewards of the actions $r_1,\dots,r_T$ with $r_t=(r_t^1,\dots,r_t^n)$ for all $t\in[T]$, where $T\in\Z_{\ge1}$, $r_t^i\in\R_{\ge0}$ is the reward of choosing $i$ in round $t$ and assume without loss of generality that $r_t^j\in[a,b]$ with $a,b\in\R$, for all $t\in[T]$ and all $i\in[n]$.  At each round $t\in[T]$, a decision maker chooses an action $i_t\in[n]$ and  receives the corresponding reward $r_t^{i_t}$. The rewards of the actions are assumed to be unknown a priori and $r_1,\dots,r_T$ are generated by an adaptive adversary, i.e., the reward $r_t^i$ can potentially depend on the previous choices $i_1,\dots,i_{t-1}$ made by the decision maker (but does not depend on $i_t$). 
For any $i_t^{\star}\in[n]$ and any $t\in[T]$, we define the following dynamic regret 
\begin{align}
R_{\M}(i^{\star})=\E\Big[\sum_{t=1}^T r_t^{i_t^{\star}}-\sum_{t=1}^T r_t^{i_t}\Big],\label{eqn:dynamic regret}
\end{align}
where the expectation is taken with respect to the randomness in the online algorithm used by the decision maker, $i_1^{\star},\dots,i_T^{\star}$ is the sequence of actions chosen by the clairvoyant, and $i^{\star}\triangleq(i^{\star}_1,\dots,i^{\star}_T)$. Similarly to \cite{auer2002nonstochastic}, we define the variability of $i_1^{\star},\dots,i_T^{\star}$ as 
\begin{equation}\label{eqn:switchings of i_t}
V_T^{i^{\star}}=1+\sum_{t=1}^{T-1}\BI\{i_t^{\star}\neq i_{t+1}^{\star}\}.
\end{equation} 

Now, we consider the scenario of the MAB problem when the decision maker observes only an erroneous version of the reward of the chosen action. Specifically, supposing an action $j\in[n]$ is chosen in any round $t\in[T]$, the decision maker observes $\beta_t(r_t^j-\epsilon_t^j+\epsilon_t)$ with an additive error $-\epsilon_t^j+\epsilon_t$ on $r_t^j$, where  $\epsilon_t^j,\epsilon_t\in\R$ and $\beta_1,\dots,\beta_T$ are i.i.d. Bernoulli random variable with parameter $\rho\in (0,1]$, i.e., $\beta_t\overset{\text{i.i.d.}}{\sim}\textbf{Bernoulli}(\rho)$. In other words, the decision maker observes the erroneous reward $r_t^j-\epsilon_t^j+\epsilon_t$ of the action $j$ chosen in round $t\in[T]$ with probability $\rho$ and observes $0$ with probability $1-\rho$.

 Based on the above discussions and notations, we introduce Algorithm~\ref{alg:exp3.S with error}, which generalizes the standard {\bf Exp3.S} algorithm \cite{auer2002nonstochastic} that considers the scenario when the decision maker observes the exact reward of the chosen action (with probability $1$), i.e., $\epsilon_t^j=\epsilon_t=0$ and $\rho=1$. Note that Algorithm~\ref{alg:exp3.S with error} chooses an action $i_t\in[n]$ in line~4 for round $t$ randomly according to the probabilities $p_t^1,\dots,p_t^n$ which are recursively updated as line~3 using the variables $\hat{r}_t^j$ and $\varpi_t^j$ given in lines~7-8 of the algorithm. We further state several remarks for Algorithm~\ref{alg:exp3.S with error} below. First, the randomness in line~4 of Algorithm~\ref{alg:exp3.S with error} is independent of the Bernoulli random variables $\beta_1,\dots,\beta_T$ described above. Next, the additive error $\epsilon_t^j$ (resp., $\epsilon_t$) is correlated with (resp., independent of) the action chosen by the algorithm in round $t\in[T]$. Finally, the settings of MAB and Algorithm~\ref{alg:exp3.S with error} described so far are all tailored to our algorithm design and regret analysis in Section~\ref{sec:OMDCO} for the OMDCO problem. In particular, we will use the following result proved in Appendix~\ref{app:proof of preliminary results} (that bounds the regret of Algorithm~\ref{alg:exp3.S with error}) to bound the regret incurred by the discrete variable in our algorithm proposed for the OMDCO problem.

\begin{prop}\label{prop:MAB with error}
Suppose there exist a constant $\underline{\epsilon}\in\R$ and $\bar{\epsilon}_t^j\in\R_{\ge0}$ such that $\underline{\epsilon}\le\epsilon_t^j\le\bar{\epsilon}_t^j$ for all $t\in[T]$ and all $j\in[n]$. Assume that $r_t^j\in[a,b]$ and $(r_t^j-\epsilon_t^j+\epsilon_t)\in[a,b]$ for all $t\in[T]$ and all $j\in[n]$, where $a,b\in\R$. Then, for any $T\ge1$ and any $\gamma\in(0,1]$, the dynamic regret of Algorithm~\ref{alg:exp3.S with error} defined in Eq.~\eqref{eqn:dynamic regret} satisfies that 
\begin{multline}\label{eqn:regret bound on exp3.S with error}
R_{\M}(i^{\star})\le(b-a)\frac{n\big(\E[V_T^{i^{\star}}]\ln(nT)+e\big)}{\gamma\rho}-\E\Big[\sum_{t=1}^T\epsilon_t^{i_t}\Big]\\-(e-2)\gamma\underline{\epsilon}T+(1-\gamma)\frac{n}{\gamma}\E\Big[\sum_{t=1}^T\bar{\epsilon}_t^{i_t}\Big]+(e-1)\gamma(b-a)T\\+(e-1)\gamma\E\Big[\sum_{t=1}^T\epsilon_t\Big],
\end{multline}
where the expectations are all taken with respect to the randomness in Algorithm~\ref{alg:exp3.S with error}.
\end{prop}

\subsection{OCO with partial information}\label{sec:OCO with bandit feedback}
Here, we study the OCO problem with partial information, including the single-point feedback and the two-point feedback described in Section~\ref{sec:problem formulation}. Similarly to our discussions in Section~\ref{sec:problem formulation}, at each round $t\in[T]$, a decision maker needs to choose a point $\bx\in\X\subseteq\R^d$ and then receives a reward given by $h_t(\bx)$, where $h_t:\R^d\to\R$. We assume that the functions $h_1(\cdot),\dots,h_T(\cdot)$ are unknown to the decision maker a priori, and the functions are generated by an adaptive adversary as we discussed in Section~\ref{sec:problem formulation}, i.e., $h_t(\cdot)$ is allowed to depend on $\bx_1,\dots,\bx_{t-1}$. We consider Algorithm~\ref{alg:bandit OCO with bandit feedback} for the OCO problem described above, where $\bg_t$ in line~7 (resp., line~11) can be viewed as an estimate of $\nabla h(\bz_t)$ based on the single-point feedback (resp., multi-point feedback).\footnote{Throughout this paper, if $\nabla h(\bx)$ does not exist for $h:\R^d\to\R$, one may replace $\nabla h(\bx)$ with a subgradient at $\bx$.} Note that Algorithm~\ref{alg:bandit OCO with bandit feedback} was first introduced in \cite{flaxman2005online,agarwal2010optimal}, but the previous results are on the static regret of the algorithm. In contrast, we prove upper bounds on the dynamic regret of Algorithm~\ref{alg:bandit OCO with bandit feedback}. Formally, for any $\bx_t^{\star}\in\X$ and any $t\in[T]$, we define the dynamic regret as
\begin{align}\label{eqn:dynamic regret of OCO}
R_{\CO}(\bx^{\star})=\E\Big[\sum_{t=1}^Th_t(\bx_t^{\star})-\sum_{t=1}^Th_t(\bx_t)\Big],
\end{align}
where the expectation is taken with respect to the randomness in Algorithm~\ref{alg:bandit OCO with bandit feedback}, $\bx_1^{\star},\dots,\bx_T^{\star}$ is the sequence of points chosen by a clairvoyant, and $\bx^{\star}\triangleq(\bx_1^{\star},\dots,\bx_T^{\star})$. Similarly to \cite{yang2016tracking}, we define the path variation of the clairvoyant's choices $\bx_1^{\star},\dots,\bx_T^{\star}$ as
\begin{equation}\label{eqn:switching of x^star}
V_T^{\bx^{\star}}=\sum_{t=1}^{T-1}\norm{\bx_{t}^{\star}-\bx_{t+1}^{\star}}.
\end{equation}
We have the following result for the dynamic regret of Algorithm~\ref{alg:bandit OCO with bandit feedback} proved in Appendix~\ref{app:proof of preliminary results}. When analyzing the regret of our algorithm proposed for OMDCO in Section~\ref{sec:OMDCO}, we use Proposition~\ref{prop:OCO with one-point evaluation} to bound the regret incurred by the continuous variable in the algorithm.
\begin{algorithm}[tb]
\textbf{Input:} Parameters $\xi\in(0,1),\delta\in\R_{>0}$, and learning rate $\eta_t\in\R_{>0}$ 
\caption{OCO with partial information}\label{alg:bandit OCO with bandit feedback}
\begin{algorithmic}[1]
\State Initialize $\bz_1\in(1-\xi)\X$ arbitrarily
\For{$t=1$ to $T$}
\State Choose $\bu_t\in\BS$ uniformly at random (u.a.r.)
\State Set $\bx_t=\bz_t+\delta\bu_t$ and $\tilde{\bx}_t=\bz_t-\delta\bu_t$
\If{single-point feedback}
\State Play $\bx_t$ and observe $h_t(\bx_t)$ 
\State Set $\bg_t=\frac{d}{\delta}h_{t}(\bx_{t})\bu_t$
\EndIf
\If{two-point feedback}
\State Play $\bx_t$ and observe $h_t(\bx_t),h_t(\tilde{\bx}_t)$
\State Set $\bg_t=\frac{d}{2\delta}(h_{t}(\bx_{t})-h_{t}(\tilde{\bx}_{t}))\bu_t$
\EndIf
\State Update $\bz_{t+1}=\Pi_{(1-\xi)\X}(\bz_t+\eta_t\bg_t)$
\EndFor
\end{algorithmic}
\end{algorithm}

\begin{prop}\label{prop:OCO with one-point evaluation}
Suppose Assumption~\ref{ass:domains of f_t} holds for $\X$, and Assumption~\ref{ass:function f_t} holds for $h_t:\R^d\to\R$ for all $t\in[T]$.\\
\noindent (a) For single point feedback, let $\eta_t=\frac{1}{T^{3/4}}$ for all $t\in[T]$, $\delta=\frac{r}{T^{1/4}}$ and $\xi=\frac{\delta}{r}$. Then, for any $T\ge2$,
\begin{multline*}
R_{\CO}(\bx^{\star})\\\le(\frac{d^2(D^2+C^2)r^2+D^2}{2}+3Gr+\frac{GD}{r}+D\E[V_T^{\bx^{\star}}])T^{3/4},
\end{multline*}
where the expectation is taken with respect to the randomness in Algorithm~\ref{alg:bandit OCO with bandit feedback}.\\
\noindent (b) For two-point feedback, let $\eta_t=\frac{1}{T^{1/2}}$ for all $t\in[T]$, $\delta=\frac{r}{T^{1/2}}$ and $\xi=\frac{\delta}{r}$. Then, for any $T\ge2$,
\begin{multline*}
R_{\CO}(\bx^{\star})\\\le(\frac{d^2(D^2+G^2)+D^2}{2}+3Gr+\frac{GD}{r}+D\E[V_T^{\bx^{\star}}])T^{1/2}.
\end{multline*}
(c) Further assume that $f_t(\cdot)$ is $\mu$-strongly concave over $\X$ for all $t\in[T]$ with $\mu\in\R_{>0}$.\footnote{A function $h:\R^d\to\R$ is $\mu$-strongly concave over $\X\subseteq\R^d$ if there exists $\mu\in\R_{\ge0}$ such that $h(\by_2)\le h(\by_1)+\nabla h(\by_1)^{\top}(\by_2-\by_1)-\frac{\mu}{2}\norm{\by_2-\by_1}^2$ for all $\by_1,\by_2\in\X$.} Consider the two-point feedback in Algorithm~\ref{alg:bandit OCO with bandit feedback}. Let $\eta_t=\frac{1}{t\mu}$ for all $t\in[T]$, $\delta=\frac{r}{T}$ and $\xi=\frac{\delta}{r}$. Then, for any $T\ge2$,
\begin{multline*}
R_{\CO}(\bx^{\star})\le (\frac{\mu G^2d^2}{2}+3Gr+\frac{GD}{r}+D\mu \E[V_T^{\bx^{\star}}] )(1+\ln T).
\end{multline*}
\end{prop}

\begin{rem}\label{remark:expected switchings}
In our OCO setup, the adaptive adversary chooses the reward function $h_t(\cdot)$ based on the decision maker's previous choices $\bx_1,\dots,\bx_{t-1}$ which are random according to Algorithm~\ref{alg:bandit OCO with bandit feedback}. Thus, $h_1(\cdot),\dots,h_T(\cdot)$ and $\bx_1^{\star},\dots,\bx_T^{\star}$ can potentially be random. Since Proposition~\ref{prop:OCO with one-point evaluation} bounds the expected regret of Algorithm~\ref{alg:bandit OCO with bandit feedback} (as per Eq.~\eqref{eqn:dynamic regret of OCO}), the regret bounds contain $\E[V_T^{\bx^{\star}}]$ with the expectation taken with respect to the randomness in $h_1(\cdot),\dots,h_T(\cdot)$ (i.e., the randomness in Algorithm~\ref{alg:bandit OCO with bandit feedback}).  Similar arguments apply to the $\E[V_T^{i^{\star}}]$ factor in the regret bounds for Algorithm~\ref{alg:exp3.S with error} provided in Proposition~\ref{prop:MAB with error} and other variation factors discussed in the next section. 
\end{rem}

\section{Algorithm Design and Regret Analysis for OMDCO}\label{sec:OMDCO}
We will split our study of the OMDCO problem into two cases in terms of the domain of the discrete variable, i.e., $H=1$ and $H>1$ in $\I=\{S\subseteq[n]:|S|\le H\}$. As we will see, our algorithm design for these two cases is different and the regret upper bound for the case $H=1$ is tighter.

\subsection{OMDCO with $\I=\{S\subseteq[n]:|S|\le1\}$}\label{sec:OMDCO with I=[n]}
We first present Algorithm~\ref{alg:OMDCO H=1} for the OMDCO problem when the domain of the discrete variable is given by $\I=\{S\subseteq[n]:|S|\le1\}$, which we may alternatively write as $\I=[n]$. In other words, for any round $t\in[T]$, the decision maker needs to choose a single element $i_t\in\I$. Note that Algorithm~\ref{alg:OMDCO H=1} includes both the single-point feedback and two-point feedback settings described in Section~\ref{sec:problem formulation}. Also note that Algorithm~\ref{alg:OMDCO H=1} leverages the {\bf Exp3.S} algorithm (Algorithm~\ref{alg:exp3.S with error}) to choose the discrete variable $i_t\in\I$ and leverages the OCO algorithm (Algorithm~\ref{alg:bandit OCO with bandit feedback}) to choose the continuous variable $\bx_t\in\X$. In particular, the {\bf Exp3.S} subroutine $\M$ in Algorithm~\ref{alg:OMDCO H=1} is applied to the instance of the MAB problem, where the set of possible actions is given by $\I$, and the reward of choosing any action $i\in\I$ in any round $t\in[T]$ is given by $f_t(i,\bx_t^{\star})$ with $\bx_t^{\star}$ defined in \eqref{eqn:optimal point}. Algorithm~\ref{alg:OMDCO H=1} feeds back $f_t(i_t,\bx_t)$ to $\M$ as the reward of choosing $i_t$ in round $t$, which can be viewed as an erroneous version of the true reward $f_t(i_t,\bx_t^{\star})$. Moreover, $\beta_t\overset{\text{i.i.d.}}{\sim}{\bf Bernoulii}(\rho)$ in the {\bf Exp3.S} subroutine $\M$ reduces to a static constant with $\rho=1$. 

We now aim to prove upper bounds on the $1$-regret of Algorithm~\ref{alg:OMDCO H=1}, which is defined in Eq.~\eqref{eqn:def of R alpha} and is denoted as $R(1)$. As we mentioned before, the major challenge in upper bounding $R(1)$ is that in the OMDCO problem, the rewards of $i_t\in\I$ and $\bx_t\in\X$ are coupled via $f_t(i_t,\bx_t)$. To tackle this challenge, we rely on the results shown in Section~\ref{sec:preliminary results} and decompose $R(1)$ into two terms, which correspond to the regret incurred by the discrete variable $i_t\in\I$ and the regret incurred by the continuous variable $\bx_t\in\X$; the proof can be found in Appendix~\ref{app:main results proofs}. 

\begin{algorithm}[tb]
\textbf{Input:} Parameters $\gamma\in(0,1],\xi\in(0,1),\delta\in\R_{>0}$, and learning rate $\eta_t\in\R_{>0}$ 
\caption{OMDCO with $\I=[n]$}\label{alg:OMDCO H=1}
\begin{algorithmic}[1]
\State Initialize an {\bf Exp3.S} subroutine $\M$ using $\gamma$
\State Initialize $\bz_1\in(1-\xi)\X$ arbitrarily
\For{$t=1$ to $T$}
\State Choose $i_t\in\I$ using $\M$
\State Choose $\bu_t\in\BS$ u.a.r.
\State Set $\bx_t=\bz_t+\delta\bu_t$ and $\tilde{\bx}_t=\bz_t-\delta\bu_t$
\If{single-point feedback}
\State Play $(i_t,\bx_t)$ and observe $f_t(i_t,\bx_t)$ 
\State Set $\bg_t=\frac{d}{\delta}f_{t}(i_t,\bx_{t})\bu_t$
\EndIf
\If{two-point feedback}
\State Play $(i_t,\bx_t)$ and observe $f_t(i_t,\bx_t),f_t(i_t,\tilde{\bx}_t)$
\State Set $\bg_t=\frac{d}{2\delta}(f_{t}(i_t,\bx_{t})-f_{t}(i_t,\tilde{\bx}_{t}))\bu_t$
\EndIf
\State Update $\bz_{t+1}=\Pi_{(1-\xi)\X}(\bz_t+\eta_t\bg_t)$
\State Feed back $f_t(i_t,\bx_t)$ to $\M$ as the reward of choosing $i_t$ in round $t$
\EndFor
\end{algorithmic}
\end{algorithm}

\begin{thm}\label{thm:regret of A_2}
Suppose that Assumptions~\ref{ass:domains of f_t}-\ref{ass:function f_t zero} hold. Consider $\hat{\bx}_t^{\star}\in\arg\max_{\bx\in\X}f_t(i_t,\bx)$ for all $t\in[T]$, where $i_t\in\I$ is chosen by Algorithm~\ref{alg:OMDCO H=1} in round $t$. Let $V_T^{i^{\star}}$ be defined as Eq.~\eqref{eqn:switchings of i_t}, where $i_t^{\star}$ is given by \eqref{eqn:optimal point} for all $t\in[T]$, and let $\hat{V}_T^{\bx^{\star}}=\sum_{t=1}^{T-1}\norm{\hat{\bx}_{t}^{\star}-\hat{\bx}_{t+1}^{\star}}$.\\
\noindent (a) Consider the single-point feedback in Algorithm~\ref{alg:OMDCO H=1}. Let $\eta_t=\frac{1}{T^{3/4}}$ for all $t\in[T]$, $\delta=\frac{r}{T^{1/4}}$, $\xi=\frac{\delta}{r}$ and $\gamma=\min\{1,\sqrt{\frac{n}{T^{1/4}}}\}$. Then, for any $T\ge2$,
\begin{multline}
R(1)\le C\sqrt{n}\big(\E[V_T^{i^{\star}}]\ln(nT)+e\big)T^{7/8}+C\sqrt{n}(2e-3)T^{7/8}\\+\sqrt{n}\big(P_1+D\E[\hat{V}_T^{\bx^{\star}}]\big)T^{7/8},
\end{multline}
where $P_1\triangleq \frac{d^2(D^2+C^2)r^2+D^2}{2}+3Gr+\frac{GD}{r}$.\\
\noindent (b) Consider the two-point feedback in Algorithm~\ref{alg:OMDCO H=1}. Let $\eta_t=\frac{1}{T^{1/2}}$ for all $t\in[T]$, $\delta=\frac{r}{T^{1/2}}$, $\xi=\frac{\delta}{r}$ and $\gamma=\min\{1,\sqrt{\frac{n}{T^{1/2}}}\}$. Then, for any $T\ge2$,
\begin{multline}
R(1)\le C\sqrt{n}\big(\E[V_T^{i^{\star}}]\ln(nT)+e\big)T^{3/4}+C\sqrt{n}(2e-3)T^{3/4}\\+\sqrt{n}\big(P_2+D\E[\hat{V}_T^{\bx^{\star}}]\big)T^{3/4},
\end{multline}
where $P_2\triangleq\frac{d^2(D^2+G^2)+D^2}{2}+3Gr+\frac{GD}{r}$.\\
\noindent (c) Further assume that $f_t(S,\cdot)$ is $\mu$-strongly concave over $\X$ for all $S\subseteq[n]$ and all $t\in[T]$ with $\mu\in\R_{>0}$. Consider the two-point feedback in Algorithm~\ref{alg:OMDCO H=1}. Let $\eta_t=\frac{1}{t\mu}$ for all $t\in[T]$,  $\delta=\frac{r}{T^{1/4}}$, $\xi=\frac{\delta}{r}$ and $\gamma=\min\{1,\sqrt{\frac{n}{T}}\}$ in Algorithm~\ref{alg:OMDCO H=1}. Then, for any $T\ge2$,
\begin{multline}
R(1)\le C\sqrt{n}\big(\E[V_T^{i^{\star}}]\ln(nT)+e\big)\sqrt{T}+C\sqrt{n}(2e-3)\sqrt{T}\\+\sqrt{n}\big(P_3+D\mu\E[\hat{V}_T^{\bx^{\star}}]\big)(1+\ln T)\sqrt{T},
\end{multline}
where $P_3\triangleq \frac{\mu G^2d^2}{2}+3Gr+\frac{GD}{r}$.
\end{thm}
$\bullet$ {\bf Choices of parameters and comparison to existing works.} Theorem~\ref{thm:regret of A_2} shows that the regret bounds on $R(1)$ are sublinear in $T$ and satisfy that $\tilde{O}(\sqrt{n}(\E[V_T^{i^{\star}}]+\E[\hat{V}^{\bx^{\star}}_T])T^{\theta})$, where $\theta\in(0,1)$, and $\tilde{O}(\cdot)$ compresses polynomial factors in $C,G,D,d,\mu,\frac{1}{r},r,\ln nT$. 
In addition, one can observe from the proof of Theorem~\ref{thm:regret of A_2} that the choice of $\gamma$ in Theorem~\ref{thm:regret of A_2} is used to balance the different terms in $T$ and $n$ in the regret bound such that the overall scaling with respect to $T$ and $n$ is optimized, e.g., the value of $\theta\in(0,1)$ is minimized. Moreover, following similar arguments to those in \cite{auer2002nonstochastic,jadbabaie2015online,yang2016tracking}, suppose there are known $B_T^i,B_T^{\bx}\in\R_{\ge0}$ such that $\E[V_{T}^{i^{\star}}]\le B_T^i$ and $\E[\hat{V}^{\bx^{\star}}_T]\le B_T^{\bx}$. Letting $\gamma=\sqrt{\frac{(B_T^i+B_T^{\bx})n}{T^{2(1-\theta)}}}$, one can show that $R(1)\le\tilde{O}(\sqrt{n(B_T^i+B_T^{\bx})}T^{\theta})$. The dynamic regret of the {\bf Exp3.S} algorithm for the MAB problem is bounded as $\tilde{O}(\sqrt{nTB_T^i})$  \cite{auer2002nonstochastic} and the dynamic regret of the OCO algorithm with two-point feedback is bounded as $O(\sqrt{TB_T^{\bx}})$ (under Assumptions~\ref{ass:domains of f_t}-\ref{ass:function f_t}) \cite{yang2016tracking}.\footnote{One can check that  $\hat{V}_T^{\bx^{\star}}$ reduces to $V_T^{\bx^{\star}}$ defined in Eq.~\eqref{eqn:switching of x^star} when there is a single continuous variable in the problem.} While the $\sqrt{n(B_T^i+B_T^{\bx})}$ factor in the regret bounds of Algorithm~\ref{alg:OMDCO H=1} for OMDCO matches with those of the existing algorithms for MAB and OCO, the regret bounds of Algorithm~\ref{alg:OMDCO H=1} (provided in Theorem~\ref{thm:regret of A_2}(a)-(b)) have worse scaling with $T$ due to the coupling of the discrete and continuous variables in the rewards.

$\bullet$ {\bf The $\E[V_T^{i^{\star}}]$ and $\E[\hat{V}_T^{\bx^{\star}}]$ factors}. The factor $\E[V_T^{i^{\star}}]$ is typically present in the dynamic regret bounds of algorithms for online optimization with a single discrete variable (e.g., \cite{auer2002nonstochastic,matsuoka2021tracking}) and measures the hardness of the clairvoyant's choices $i_1^{\star},\dots,i_T^{\star}$ that is a consequence of the power of the adversary who can choose the function sequence $f_1(\cdot,\cdot),\dots,f_T(\cdot,\cdot)$ in an adaptive manner (see our discussions in Section~\ref{sec:problem formulation}). 
The $\E[\hat{V}^{\bx^{\star}}_T]$ factor depends on $f_1(\cdot,\cdot),\dots,f_T(\cdot,\cdot)$ and the {\bf Exp3.S} subroutine $\M$'s choices $i_1,\dots,i_T$ in Algorithm~\ref{alg:OMDCO H=1}. Suppose that the adversary switches the function $f_t(\cdot,\cdot)$ for at most $\lambda\in\Z_{\ge1}$ times, and the switching points $t=t_1,\dots,t_{\lambda}$ satisfy $t_j-t_{j-1}=T/\lambda$ for all $j\in[\lambda+1]$ with $t_0=1$ and $t_{\lambda+1}=T$.\footnote{For simplicity, we assume that $T/\lambda\in\Z$; otherwise we need to consider $t_{j+1}-t_j=\lceil T/\lambda \rceil$ for all $j\in[\lambda]$ and $t_{\lambda+1}=T$.} Meanwhile, we let the {\bf Exp3.S} subroutine $\M$ in Algorithm~\ref{alg:OMDCO H=1} only switch its choice $i_t$ at rounds $t=t_1,\dots,t_{\lambda}$, which implies that $\E[\hat{V}^{\bx^{\star}}_T]\le 2(\lambda+1) D$. Considering the setting in Theorem~\ref{thm:regret of A_2}(c), and following similar arguments to those for  \cite[Theorem~14]{altschuler2018online} and \eqref{eqn:R_Ae inter} in the proof of Theorem~\ref{thm:regret of A_2}, one can show that
\begin{multline*}
R(1)\le \frac{CT}{\lambda}\frac{n\big(\E[V_{T^{\lambda}}^{i^{\star}}]\ln(n\lambda)+e\big)}{\gamma}\\+\frac{n}{\gamma}\big(P_3+2\mu D^2(\lambda+1)\big)(1+\ln \lambda)+C(2e-3)\gamma T,
\end{multline*}
where $V_{T^{\lambda}}^{i^{\star}}=1+\sum_{t=t_1}^{t_{\lambda}}\BI\{i_t^{\star}\neq i_{t+1}^{\star}\}$. Letting $\lambda=\lceil T^{1/2}\rceil$ and  $\gamma=\sqrt{\frac{n}{T^{1/2}}}$, we obtain $R(1)\le \tilde{O}(\sqrt{n}\E[V_{T^{\lambda}}^{i^{\star}}]T^{3/4})$.

\subsection{OMDCO with $\I=\{S\subseteq[n]:|S|\le H\}$}\label{sec:OMDCO with general I}
Now, we present Algorithm~\ref{alg:OMDCO H >1} for the OMDCO problem when the domain of the discrete variable is given by the uniform matroid $\I=\{S\subseteq[n]:|S|\le H\}$ with $H\in\Z_{\ge2}$. Similarly to Algorithm~\ref{alg:OMDCO H=1}, Algorithm~\ref{alg:OMDCO H >1} leverages the {\bf Exp3.S} algorithm (Algorithm~\ref{alg:exp3.S with error}) and the OCO algorithm (Algorithm~\ref{alg:bandit OCO with bandit feedback}) to choose the discrete variable $S_t\in\I$ and the continuous variable $\bx_t\in\X$, respectively. Different from Algorithm~\ref{alg:OMDCO H=1}, Algorithm~\ref{alg:OMDCO H >1} deals with the case when the domain of the discrete variable is a collection of subsets of $[n]$ that has cardinality no greater than $H$. Thus, Algorithm~\ref{alg:OMDCO H >1} initializes $H$ independent copies of the {\bf Exp3.S} algorithm $\M_1,\dots,\M_H$, let each $\M_l$ ($l\in[H]$) choose a single element $s_t^l$ in round $t\in[T]$, and constructs the set $S^{\prime}_t=\{s_t^l:l\in[H]\}$. Denote $S^{\prime l}_t=\{s^{\prime 1},\dots,s^{\prime l}\}$ for all $l\in[H]$ with $S^{\prime 0}=\emptyset$. For any $l\in[H]$, the {\bf Exp3.S} subroutine $\M_l$ in Algorithm~\ref{alg:OMDCO H >1} is applied to the instance of MAB, where the set of possible actions is given by $[n]$, and the reward of any action $i\in[n]$ in round $t\in[T]$ is given by $f_t(S_{t}^{\prime l-1}\cup\{i\},\bx^{\star}_t)-f_t(S^{\prime l-1}_t,\bx_t^{\star})$. We then see from line~12 of Algorithm~\ref{alg:OMDCO H >1} that the algorithm feeds back $\BI\{\tilde{\beta}_t=1,l_t=l,i_t=s_t^l\}f_t(S_t,\bx_t)$ to $\M_l$ as the reward of choosing $s_t^l\in[n]$ in round $t\in[T]$, where $f_t(S_t,\bx_t)$ can be viewed as an erroneous version of the true reward $f_t(S_{t}^{\prime l-1}\cup\{s_t^l\},\bx^{\star}_t)-f_t(S^{\prime l-1}_t,\bx_t^{\star})$. 

\begin{algorithm}[tb]
\textbf{Input:} Parameters $\gamma,\tilde{\rho}\in(0,1],\xi\in(0,1),\delta\in\R_{>0}$, and learning rate $\eta_t\in\R_{>0}$ 
\caption{OMDCO with $\I=\{S\subseteq[n]:|S|\le H\}$}\label{alg:OMDCO H >1}
\begin{algorithmic}[1]
\State Initialize $H$ independent copies of the {\bf Exp3.S} algorithm, denoted as $\M_1,\dots,\M_H$,  using $\gamma$. For any $l\in[H]$, use $\M_l$ to choose $s_1^l\in[n]$ and set $S^{\prime}_1=\{s_1^l:l\in[H]\}$
\State Initialize $\bz_1\in(1-\xi)\X$ arbitrarily
\For{$t=1$ to $T$}
\State Sample $\tilde{\beta}_t\overset{\text{i.i.d.}}{\sim}{\bf Bernoulli}(\tilde{\rho})$
\State Choose $l_t\in[H]$ u.a.r. and choose $i_t\in[n]$ u.a.r. 
\State Choose $\bu_t\in\BS$ u.a.r.
\State Set $S_t=\begin{cases}S^{\prime l_t-1}_{t}\cup\{i_t\},\ \text{if}\ \tilde{\beta}_t=1\\ S^{\prime}_t,\ \text{if}\ \tilde{\beta}_t=0\end{cases}$
\State Set $\bx_t=\bz_t+\delta\bu_t$ and $\tilde{\bx}_t=\bz_t-\delta\bu_t$
\State Play $(S_t,\bx_t)$ and  observe $f_t(S_t,\bx_t),f_t(S_t,\tilde{\bx}_t)$
\State Set $\bg_t=\frac{d}{2\delta}(f_{t}(S_t,\bx_{t})-f_{t}(S_t,\tilde{\bx}_{t}))\bu_t$
\State Update $\bz_{t+1}=\Pi_{(1-\xi)\X}(\bz_t+\eta_t\bg_t)$
\State For any $l\in[H]$, feed back $\BI\{\tilde{\beta}_t=1,l_t=l,i_t=s_t^l\}f_t(S_t,\bx_t)$ to $\M_l$ as the reward of choosing $s_t^l$ in round $t$
\State For any $l\in[H]$, use $\M_l$ to choose $s_{t+1}^l$ and set $S^{\prime}_{t+1}=\{s_{t+1}^l:l\in[n]\}$
\EndFor
\end{algorithmic}
\end{algorithm}

Recall that we prove upper bounds on the $1$-regret of Algorithm~\ref{alg:OMDCO H=1} in Section~\ref{sec:OMDCO with I=[n]}, since Algorithm~\ref{alg:OMDCO H=1} deals with the instances of the OMDCO problem when the discrete domain $\I=[n]$ has no combinatorial structure. In contrast, Algorithm~\ref{alg:OMDCO H >1} deals with the OMDCO instances when the discrete domain $\I=\{S\subseteq[n]:|S|\le H\}$ with $H\in\Z_{>1}$ has a combinatorial nature. In fact, the offline problem of $\max_{S\subseteq[n],|S|\le H} g(S)$ (with a set function $g:2^{[n]}\to\R$ and $H\in\Z_{>1}$) is NP-hard in general (see e.g. \cite{feige1998threshold}), i.e., any polynomial-time algorithm only returns an approximately optimal solution (unless P$=$NP). Thus, when characterizing the performance of Algorithm~\ref{alg:OMDCO H >1} for the OMDCO problem with $\I=\{S\subseteq[n]:|S|\le H\}$, we consider the $\alpha$-regret of the algorithm defined in Eq.~\eqref{eqn:def of R alpha} with $\alpha\in(0,1]$, which compares the reward of the decisions $(S_1,\bx_1),\dots,(S_T,\bx_T)$ against $\alpha$ times the reward of an optimal decision sequence $(S_1^{\star},\bx_1^{\star}),\dots,(S_T^{\star},\bx_T^{\star})$ defined in~\eqref{eqn:optimal point}. Note that the notion of $\alpha$-regret has also been used to characterize the performance of algorithms for online optimization with a single discrete variable \cite{streeter2008online}. In the sequel, we denote the $\alpha$-regret of Algorithm~\ref{alg:OMDCO H >1} as $R_{A_u}(\alpha)$. Note that the value of $\alpha$ depends on certain parameters of the functions $f_t(\cdot,\cdot)$. To this end, we introduce the following definitions (see, e.g., \cite{bian2017guarantees,das2018approximate}),

\begin{defn}
\label{def:submodularity ratio}
The submodularity ratio of a set function $g:2^{[n]}\to\R_{\ge0}$ is the largest $\kappa_g\in\R$ such that 
\begin{equation*}
\sum_{\omega\in \Omega\setminus S}\big(g(S\cup\{\omega\})-g(S)\big)\ge\kappa_g\big(g(S\cup \Omega)-g(S)\big),
\end{equation*}
for all $S,\Omega\subseteq[n]$.
\end{defn}
\begin{defn}
\label{def:curvature}
The curvature of a set function $g:2^{[n]}\to\R_{\ge0}$ is the smallest $c_g\in\R$ such that
\begin{equation*}
g(\Omega\cup\{\omega\})-g(\Omega)\ge(1-c_g)\big(g(S\cup\{\omega\})-g(S)\big),
\end{equation*} 
for all $S\subseteq\Omega\subseteq[n]$ and all $\omega\in[n]\setminus\Omega$.
\end{defn}
For monotone nondecreasing $g(\cdot)$, one can check that $\kappa_g,c_g\in[0,1]$ for all $t\in[T]$. If we assume that $g(\cdot)$ is also submodular,\footnote{A set function $h:2^{[n]}\to\R$ is submodular if and only if $\sum_{\omega\in\Omega\setminus S}\big(h(S\cup\{\omega\})-h(S)\big)\ge h(S\cup\Omega)-h(S)$} one can show via Definition~\ref{def:submodularity ratio} that $\gamma_g=1$ (e.g., \cite{bian2017guarantees,das2018approximate}). If we assume that $g(\cdot)$ is modular,\footnote{A set function $g:2^{[n]}\to\R_{\ge0}$ is modular if and only if $g(\Omega)=\sum_{\omega\in\Omega}g(\omega)$ for all $\Omega\subseteq[n]$.} we see from Definition~\ref{def:curvature} that $c_g=0$. In words, the submodularity ratio of a monotone nondecreasing set function characterizes its approximate submodularity; the curvature of a set function characterizes how far the function $g(\cdot)$ is from being modular. We will further make the following assumption, which says that $\bx_t^{\star}$ (from an optimal decision point $(S_t^{\star},\bx_t^{\star})$) is optimal over all $S\subseteq[n]$ and all $\bx\in\X$. Again, we will provide examples in Section~\ref{sec:applications} that satisfy this assumption.
\begin{assum}\label{ass:x_t^star optimal}
For any $t\in[T]$ and any $S\subseteq[n]$, there exists $\bx_t^{\star}$ given by \eqref{eqn:optimal point} such that $f_t(S,\bx_t^{\star})\ge f_t(S,\bx)$ for all $\bx\in\X$.
\end{assum}
Similarly to Eq.~\eqref{eqn:switchings of i_t}, we define the variability of $S_1^{\star},\dots,S_T^{\star}$ as 
\begin{equation}\label{eqn:switchings of S_t}
V_T^{S^{\star}}=1+\sum_{t=1}^{T-1}\BI\{S_t^{\star}\neq S_{t+1}^{\star}\}.
\end{equation}
We now prove the following upper bounds on the regret of Algorithm~\ref{alg:OMDCO H >1}; the proof can be found in Appendix~\ref{app:main results proofs}.
\begin{thm}\label{thm:regret of A_3}
Suppose that Assumptions~\ref{ass:domains of f_t}-\ref{ass:function f_t zero} and \ref{ass:x_t^star optimal} hold. Let $x_t^{\star}\in\X$ given by \eqref{eqn:optimal point} satisfy Assumption~\ref{ass:x_t^star optimal} for all $t\in[T]$. Consider $\alpha=\frac{1}{\overline{c}}(1-e^{-\overline{c}\underline{\kappa}})$, where $\underline{\kappa}\triangleq\min_{t\in[T]}\kappa_t$ and $\overline{c}\triangleq\max_{t\in[T]}c_t$ with $\kappa_t$ and $c_t$ to be the submodularity ratio and curvature of $f_t(\cdot,\bx^{\star}_t)$ given by Definitions~\ref{def:submodularity ratio} and~\ref{def:curvature}, respectively.\\
\noindent(a) Let $\eta_t=\frac{1}{T^{1/2}}$ for all $t\in[T]$,  $\delta=\frac{r}{T^{1/2}}$, $\xi=\frac{\delta}{r}$, $\gamma=\min\{1,\frac{n}{T^{1/6}}\}$ and $\tilde{\rho}=\min\{1,\sqrt{\frac{H}{T^{1/3}}}\}$ in Algorithm~\ref{alg:OMDCO H >1}. Then, for any $T\ge2$,
\begin{multline*}
R(\alpha)\le C\sqrt{H}\big(\E[V_T^{S^{\star}}]\ln(nT)+e\big)T^{1/3}+C\sqrt{H}T^{5/6}\\+Cn(3e-4)T^{5/6}+\sqrt{H}n(P_2+D\E[V_T^{\bx^{\star}}])T^{5/6},
\end{multline*}
where $P_2=\frac{d^2(D^2+G^2)+D^2}{2}+3Gr+\frac{GD}{r}$ and $V_T^{\bx^{\star}}$ is defined as Eq.~\eqref{eqn:switching of x^star}. \\
\noindent(b) Further assume that $f_t(\cdot)$ is $\mu$-strongly concave for all $t\in[T]$ with $\mu\in\R_{>0}$. Let $\eta_t=\frac{1}{t\mu}$ for all $t\in[T]$, $\delta=\frac{r\ln T}{T}$, $\xi=\frac{\delta}{r}$, $\gamma=\min\{1,\frac{n}{T^{1/3}}\}$ and $\tilde{\rho}=\min\{1,\frac{\sqrt{H}}{T^{1/3}}\}$ in Algorithm~\ref{alg:OMDCO H >1}. Then, for any $T\ge2$,
\begin{multline*}
R(\alpha)\le C\sqrt{H}\big(\E[V_T^{S^{\star}}]\ln(nT)+e\big)T^{2/3}+C\sqrt{H}T^{2/3}\\+Cn(3e-4)T^{2/3}+2\sqrt{H}n(P_3+D\mu\E[V_T^{\bx^{\star}}])(1+\ln T)T^{2/3},
\end{multline*}
where $P_3=\frac{\mu G^2d^2}{2}+3Gr+\frac{GD}{r}$.
\end{thm}

$\bullet$ {\bf Choices of parameters and factors in the regret bounds.}
The bounds on $R(\alpha)$ provided in Theorem~\ref{thm:regret of A_3} are sublinear in $T$ and satisfy that $\tilde{O}(n\sqrt{H}(\E[V_T^{S^{\star}}]+\E[V_T^{\bx^{\star}}])T^{\theta})$ for some $\theta\in(0,1)$. 
Similarly to our arguments for Theorem~\ref{thm:regret of A_2} in Section~\ref{sec:OMDCO with I=[n]}, the choices of $\gamma$ and $\rho$ in Theorem~\ref{thm:regret of A_3} are used to balance the different terms in $H,n,T$ in the regret bound, which yields the optimal scaling with respect to $H,n,T$ in the overall regret bound. Moreover, letting $\gamma=n\sqrt{\frac{(B_T^S+B_T^{\bx})}{T^{2(1-\theta)}}}$, one can show that $ R(\alpha)\le\tilde{O}(n\sqrt{H(B_T^S+B_T^{\bx})}T^{\theta})$, where $B_T^S,B_T^{\bx}\in\R_{\ge0}$ are known upper bounds with  $\E[V_{T}^{S^{\star}}]\le B_T^{S}$ and $\E[V_T^{\bx^{\star}}]\le B_T^{\bx}$. The term $E[V_T^{S^{\star}}]$ (resp., $\E[V_T^{\bx^{\star}}]$) measures the hardness of the clairvoyant's choice $S_1^{\star},\dots,S_T^{\star}$ (resp., $\bx_1^{\star},\dots,\bx_T^{\star}$), and can be sublinear in $T$ if the clairvoyant switches its choice $S_t^{\star}$ (resp., $\bx_t^{\star}$) for a limited number of times. Finally, if the set function $f_t(\cdot,\bx_t^{\star})$ is modular, i.e., $c_t=0$, for all $t\in[T]$, then $\overline{c}=0$, which implies $\alpha=1$. 

$\bullet$ {\bf Comparison to existing works.} 
Under the bandit information setting (i.e., single-point feedback), \cite{streeter2008online} gives an online algorithm for maximization of a submodular set function subject to a cardinality constraint and upper bounds its {\it static} $(1-1/e)$-regret as $O(H(n\ln n)^{1/3}T^{2/3})$. When considering a partition matroid constraint,\footnote{Given a partition of ground set $\V=\V_1\cap\cdots\cap\V_H$ with $|\V|=n$, a partition matroid is a collection $\mathcal{C}$ of subsets of $\V$ s.t. for any $S\in\mathcal{C}$, $|S\cap\V_k|\le 1$ $\forall k\in[H]$.} \cite{golovin2014online} provides an online algorithm with $O(H^{3/2}(n^2\ln n)^{1/3}T^{3/4})$ static $(1-1/e)$-regret. Thus, the scaling with $n$ and $H$ in the regret bounds of Algorithm~\ref{alg:OMDCO H >1} provided in Theorem~\ref{thm:regret of A_3} is comparable to that in the regret bounds of the algorithms described above. When allowing two-point feedback, we show in Theorem~\ref{thm:regret of A_3}(b) that the regret bound of Algorithm~\ref{alg:OMDCO H >1} achieves the same scaling with $T$ as the algorithm in \cite{streeter2008online}. Further, our algorithm can handle nonsubmodularity in the reward function $f_t(\cdot,\cdot)$ via Definitions~\ref{def:submodularity ratio}-\ref{def:curvature}.

$\bullet$ {\bf Relaxation of Assumption~\ref{ass:x_t^star optimal}.} We argue that Assumption~\ref{ass:x_t^star optimal} can be relaxed to $\sum_{t=1}^T\big(f_t(S_t,\bx_t^{\star})-f_t(S_t,\bx_t)\big)\ge-\tilde{O}(\E[V_T^{\bx^{\star}}]T^{\theta})$, where $(S_t,\bx_t)$ is chosen by Algorithm~\ref{alg:OMDCO H >1} for all $t\in[T]$. In fact, supposing the function $f_t(S_t,\bx)$ is linear with respect to $\bx$, which implies that $f_t(S_t,\bx)$ is also convex with respect to $\bx$, one can follow similar arguments to those for Proposition~\ref{prop:OCO with one-point evaluation} and show that the relaxed assumption described above holds. Based on the relaxed assumption, one can then use similar arguments to those for Theorem~\ref{thm:regret of A_3} and show that the upper bounds on $ R(\alpha)$ will become $\tilde{O}(\E[V_T^{\bx^{\star}}+\hat{V}_T^{\bx^{\star}}+V_T^{S^{\star}}]T^{\theta})$, where $\E[\hat{V}_T^{\bx^{\star}}]=\sum_{t=1}^{T-1}\norm{\hat{\bx}_t^{\star}-\hat{\bx}_{t+1}^{\star}}$ with $\hat{\bx}^{\star}_t\in\argmax_{\bx\in\X}f_t(S_t,\bx)$ for all $t\in[T]$, and $S_1,\dots,S_T$ are chosen by the {\bf Exp3.S} subroutines $S_1,\dots,S_T$ in Algorithm~\ref{alg:OMDCO H >1}.

\section{Applications}\label{sec:applications}
As an application, we focus on the subset selection problem, whose offline setting has been widely studied in the literature (e.g., \cite{das2008algorithms,wei2015submodularity,elenberg2018restricted}) and is given by the following general form: 
\begin{equation}\label{eqn:obj ss}
	\max_{\substack{\bx\in\X,\text{supp}(\bx)\le H}}f^s(\bx),
\end{equation}
where $\X\subseteq\R^n$, $H\in\Z_{\ge1}$, $f^s:\R^n\to\R_{\ge0}$ with $f^s(\mathbf{0})=0$. Introducing a discrete variable $S\subseteq[n]$, \eqref{eqn:obj ss} may be equivalently written as
\begin{equation}
	\max_{\substack{S\subseteq[n],|S|\le H\\ \bx\in\X}}f^s(\bx(S)),
\end{equation}
Hence, we may also equivalently view $f^s (\cdot)$ as $f^s(\cdot,\cdot)$, where $f^s:2^{[n]}\times\R^n\to\R_{\ge0}$. The subset selection problem was initially motivated by sparse modeling in data analysis, where $[n]$ is the set of all features (i.e., attributes) of the data points and the problem is to build a model using a subset of features of size at most $H$. Many other real-world applications fit into the framework of the subset selection problem. For instance, in influence maximization over social networks (e.g., \cite{alon2012optimizing,miyauchi2015threshold}), we need to choose a subset of initially activated nodes in the network and determine the influence level of each activated node such that the impact on the remaining nodes is maximized. For network protection against a virus spreading process (e.g., \cite{preciado2014optimal}), one needs to choose a subset of the nodes and determine the amount of vaccinations and antidotes for each chosen node to minimize the spread of the virus over the network. In the problem of sensor selection for state or parameter estimation, one needs to select a subset of all the candidate sensors to collect measurements and determine the resource (e.g., sensing or communication power) allocated to each selected sensor such that the estimation performance is optimized (e.g., \cite{mo2011sensor,ye2020complexity}). 

The works mentioned above consider the {\it offline} setting of problem~\eqref{eqn:obj ss}, i.e., the objective function $f^s(\cdot)$ is known a priori and does not change over time. However, in many real-world applications (such as the ones described above), the underlying environment is typically unknown to the decision maker and may also change over time (e.g., the virus spreading process with an unknown and time-varying infection rate). If the decision maker does not know the objective function and needs to interact with the unknown and time-varying environment while making decisions, the Online Subset Selection (OSS) problem results. Specifically, at each round $t\in[T]$, $S_t\subseteq[n]$ (with $|S_t|\le H$) denotes the set of chosen elements in round $t\in[T]$, and $f^s_t(\bx_t(S_t))$ (defined similarly to $f^s(\cdot)$ above) returns the corresponding reward and is assumed to be unknown (see our discussions in Section~\ref{sec:problem formulation}). Thus, the OSS problem is a special class of the OMDCO problem. As we argued in the previous sections, one can find instances of the OSS problem that satisfy Assumptions~\ref{ass:function f_t zero} and \ref{ass:x_t^star optimal}. To this end, we first prove the following result.



\begin{lem}\label{lemma:extra assumption for OSS}
		Consider any $S\subseteq\Omega\subseteq[n]$ and any $t\in[T]$, and let $\bx_t^{\star\Omega}\in\argmax_{\bx\in\X}f_t^s(\Omega,\bx)$. Suppose $f_t^s(\Omega,\bx_t^{\star\Omega})\ge f_t^s(S,\bx_t^{\star\Omega})\ge f_t^s(S,\bx_t^{\star S})$. Then, Assumptions~\ref{ass:function f_t zero} and \ref{ass:x_t^star optimal} hold for $f_t^s(\cdot,\cdot)$.
	\end{lem}
	
	Observing that $\bx_t^{\star}=\bx_t^{\star[n]}$ satisfies \eqref{eqn:optimal point}, the proof of Lemma~\ref{lemma:extra assumption for OSS} follows by verification. We now give the following examples for the OSS problem described above (with the corresponding objective function  $f_t^s(\cdot,\cdot)$). Using Lemma~\ref{lemma:extra assumption for OSS}, one can verify that $f_t^s(\cdot,\cdot)$ in Examples~\ref{exp:OSS 1}-\ref{exp:OSS 2} satisfy Assumptions~\ref{ass:function f_t zero} and \ref{ass:x_t^star optimal}.

\begin{exmp}\label{exp:OSS 1}
	For any $S\subseteq[n]$, any $\bx\in\R^n$ and any $t\in[T]$, define $f_t^s (S,\bx)=h_t^{s}(\sum_{i\in S}h^i_t((\bx)_i)$, where $h_t^s:\R_{\ge0}\to\R_{\ge0}$ is monotone nondecreasing and $h_t^i:\R\to\R_{\ge0}$ for all $i\in[n]$. Let the continuous domain satisfy that $\X=\X_1\times\cdots\times\X_n$, i.e., $(\bx)_i\in\X_i$ for all $i\in[n]$.
\end{exmp}

\begin{exmp}\label{exp:fcl}
	For any $S\subseteq[n]$, any $\bx\in\R^n$ and any $t\in[T]$, define $f_t^s(S,\bx)=\sum_{i=1}^n\sum_{j\in S}h_t^{i,j}((\bx)_j)$, where $h_t^{i,j}:\R\to\R_{\ge0}$. Let the continuous domain satisfy that $\X=\X_1\times\cdots\times\X_n$.
\end{exmp}

\begin{exmp}\label{exp:OSS 2}
	For any $S\subseteq[n]$, any $\bx\in\R^n$ and any $t\in[T]$, define $f_t^s(S,\bx)=\bw_t^{\top}\bx(S)+\mathbf{b}_t$, where $\bw_t,\mathbf{b}_t\in\R^n$. Let $\bw_t$ satisfy that $(\bw_t)_i\ge0$ for all $i\in[n]$ and let $\X$ satisfy that $(\bx)_i\ge0$ for all $\bx\in\X$.
\end{exmp}

Example~\ref{exp:fcl} can be translated to the facility location problem whose offline and online settings have been widely studied in the literature (see, e.g., \cite{meyerson2001online,drezner2004facility,adibi2022minimax}, for more details about the problem setup). Also note that Examples~\ref{exp:OSS 1}-\ref{exp:fcl} require $\X=\X_1\times\cdots\times\X_n$ to have the disjoint structure such that Assumptions~\ref{ass:function f_t zero} and \ref{ass:x_t^star optimal} are satisfied. To consider more general structure of the continuous domain $\X$, we recall from our discussions for Theorem~\ref{thm:regret of A_3} in Section~\ref{sec:OMDCO with general I} that Assumption~\ref{ass:x_t^star optimal} can be relaxed to $\sum_{t=1}^T\big(f_1(S_t,\bx_t^{\star})-f_t(S_t,\bx_t)\big)\ge-\tilde{O}(\E[V_T^{\bx^{\star}}]T^{\theta})$ for some $\theta\in(0,1)$, which is satisfied by $f_t(S_t,\bx)$ that is linear with respect to $\bx$. Under the linear setting, one can verity that $f_t^s(\cdot,\cdot)$ with $\X$ described in Example~\ref{exp:OSS 2} also satisfy Assumption~\ref{ass:function f_t zero}. The linear function model considered in Example~\ref{exp:OSS 2} is an important special case in online optimization problems (e.g., \cite{bubeck2011introduction,abernethy2009beating,li2022online}). Further supposing $f_t^s(S,\bx)=\sum_{i\in S}h_t^i((\bx)_i)$ in Example~\ref{exp:OSS 1}, one can show that $f_t^s(\cdot,\bx_t^{\star})$ is a modular set function. Similarly, one can check that $f_t^s(\cdot,\bx_t^{\star})$ given in Examples~\ref{exp:fcl}-\ref{exp:OSS 2} are also modular.  Since Examples~\ref{exp:OSS 1}-\ref{exp:OSS 2} are instances of the OSS problem, they can be used to model the different real-world applications discussed before.

Supposing Assumptions~\ref{ass:domains of f_t}-\ref{ass:function f_t} also hold, one can now apply Algorithm~\ref{alg:OMDCO H=1} to the OSS problem with $H=1$ and achieve the guarantee on the $1$-regret given by Theorem~\ref{thm:regret of A_2}. To apply Algorithm~\ref{alg:OMDCO H >1} to OSS with $H\in\Z_{>1}$ and achieve the guarantee on the $\alpha$-regret given in Theorem~\ref{thm:regret of A_3}, we may leverage the following results, which work for different classes of the function $f_t^s(\cdot,\cdot)$. The proof of Proposition~\ref{prop:sub ratio and curvature in SS} follows from similar arguments to those for  \cite[Theorem~1]{elenberg2018restricted} and is omitted. The proof of Proposition~\ref{prop:sub ratio and curvature in SS 2} can be found in Appendix~\ref{app:application proof}.
\begin{prop}\label{prop:sub ratio and curvature in SS}
	Consider the setting in Lemma~\ref{lemma:extra assumption for OSS} and suppose $f_t^s(S,\cdot):\R^n\to\R_{\ge0}$ is $\sigma_t$-smooth and $\mu_t$-strongly concave over the domain $\X$ for all $S\in\I$,\footnote{A function $h:\R^n\to\R$ is $\sigma$-smooth over $\X\subseteq\R^n$ if there exists $\sigma_t\in\R_{\ge0}$ such that $h(\by_2)\ge h(\by_1)+\nabla h(\by_1)^{\top}(\by_2-\by_1)-\frac{\sigma}{2}\norm{\by_2-\by_1}^2$ for all $\by_1,\by_2\in\X$.} with $\sigma_t,\mu_t\in\R_{>0}$, for all $t\in[T]$. 
	Then, the submodularity ratio and curvature of $f_t^s(\cdot,\bx_t^{\star})$, given by Definitions~\ref{def:submodularity ratio} and \ref{def:curvature} and denoted as $\kappa_t^s$ and $c_t^s$, satisfy that $\kappa_t^s\ge\underline{\kappa}_t^s$ and $c_t^s\le 1-\underline{\kappa}_t^s$, respectively, where $\bx_t^{\star}=\bx_t^{[n]_{\star}}$ with $\bx_t^{[n]_{\star}}$ given in Assumption~\ref{ass:x_t^star optimal} and $\underline{\kappa}_t^s=\mu_t/\sigma_t$.
\end{prop}

\begin{prop}\label{prop:sub ratio and curvature in SS 2}
	Suppose  Assumption~\ref{ass:function f_t zero} holds and $f_t^s(S,\bx)$ is differentiable and concave with respect to $\bx$ for all $S\subseteq[n]$ and all $t\in[T]$. Then, the submodularity ratio of $f_t^s(\cdot,\bx_t^{\star})$ satisfies that $\kappa_t^s\ge\underline{\kappa}_t^s$, where $\bx_t^{\star}$ is defined as\eqref{eqn:optimal point} and $\underline{\kappa}_t^s=\min_{\omega\in[n]}\frac{\min_{\bx\in\X}|(\nabla f_t^s(\bx))_{\omega}|}{\max_{\bx\in\X}|(\nabla f_t^s(\bx))_{\omega}|}$.
\end{prop}

\begin{rem}
Further assume that $f_t^s(\bx)$ is linear with respect to $\bx$, we see that $\underline{\kappa}_t^s=1$ in Proposition~\ref{prop:sub ratio and curvature in SS 2}. 
\end{rem}
Recall from Theorem~\ref{thm:regret of A_3} that $\alpha=\frac{1}{\overline{c}}(1-e^{-\overline{c}\underline{\kappa}})$. Under Assumptions~\ref{ass:domains of f_t}-\ref{ass:function f_t} and considering the setting in Proposition~\ref{prop:sub ratio and curvature in SS}, one can show that $\alpha^s=\frac{1}{1-\underline{\kappa}^s}(1-e^{-(1-\underline{\kappa}^s)\underline{\kappa}^s})\le\alpha$, where $\underline{\kappa}^s=\min_{t\in[T]}\underline{\kappa}_t^s$, and apply Algorithm~\ref{alg:OMDCO H >1} to the OSS problem with $H\in\Z_{>1}$ to achieve the $\alpha^s$-regret upper bounded by Theorem~\ref{thm:regret of A_3}. Similarly, supposing Assumptions~\ref{ass:domains of f_t}-\ref{ass:function f_t zero} and \ref{ass:x_t^star optimal} hold and considering the setting in Proposition~\ref{prop:sub ratio and curvature in SS 2}, one can show that $\alpha^s=1-e^{-\underline{\kappa}^s}\le\alpha$, where we use the naive upper bound $c_t^s\le 1$.

\section{Numerical Results}
We validate our theoretical results using numerical experiments using Example~\ref{exp:OSS 1} with $\I=\{S\subseteq[n]:|S|\le H\}$. Let $n=5$, $H=3$ and $\X_i=[-1,4]$ for all $i\in[n]$. For any $t\in[T]$, let $f_t^s(S,\bx)=\sum_{i\in S}(-a_t^i(\bx)_i^2+b_t^i(\bx)_i+c_t^i)$ for all $S\subseteq[n]$ and all $\bx\in\X$, where we draw $a_t^i,b_t^i\in[1,4]$ randomly with $c_t^i=70$ fixed over all $t\in[T]$ and all $i\in[n]$. As we argued in Section~\ref{sec:applications}, $f_t^s(\cdot,\bx_t^{\star})$ is modular, which implies that $\alpha$ considered in Theorem~\ref{thm:regret of A_3} satisfies $\alpha=1$. Moreover, one can check that $f_t^s(S,\cdot)$ is $2$-strongly concave for all $S\in\I$ with $|S|=H$. 

\begin{figure}[htbp]
	\centering
	\subfloat[a][]{\includegraphics[width=0.5\linewidth]{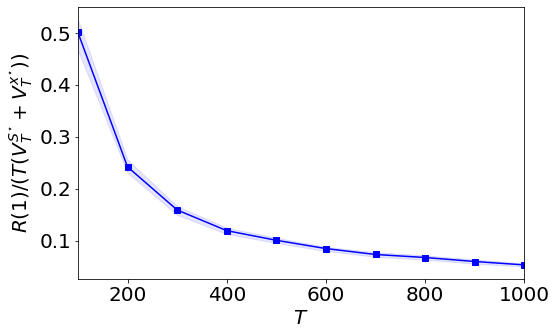}} 
	\subfloat[b][]{\includegraphics[width=0.5\linewidth]{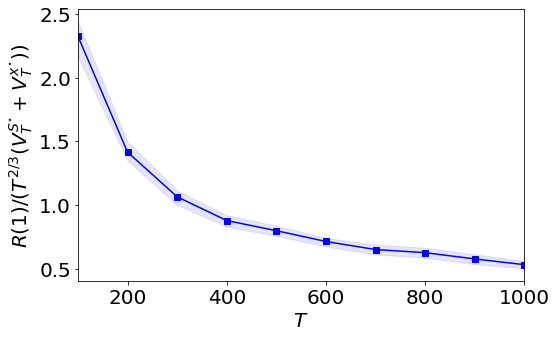}}
	\caption{Performance of Algorithm~\ref{alg:OMDCO H >1} when applied to the constructed OSS instances. Shaded regions represent quantiles.}
	\label{fig:dynamic regret}
\end{figure}

\begin{figure}[htbp]
	\centering
	\subfloat[a][]{\includegraphics[width=0.5\linewidth]{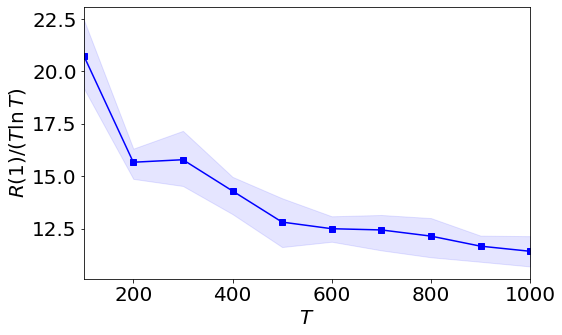}} 
	\subfloat[b][]{\includegraphics[width=0.5\linewidth]{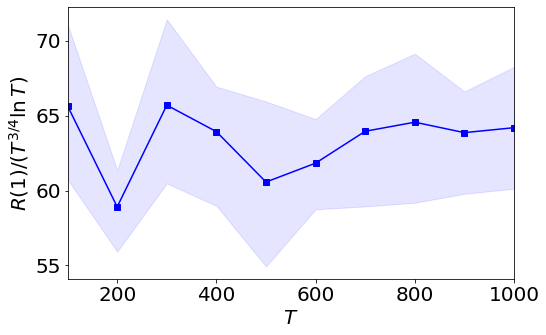}}
	\caption{Performance of Algorithm~\ref{alg:OMDCO H >1} when applied to the constructed OSS instances with limited switchings in $f_1^s(\cdot,\cdot),\dots,f_T^s(\cdot,\cdot)$. Shaded regions represent quantiles.}
	\label{fig:dynamic regret limited}
\end{figure}

We apply Algorithm~\ref{alg:OMDCO H >1} to the OSS instance constructed above; all the results are averaged over $50$ experiments. In Fig.~\ref{fig:dynamic regret}, we plot $ R(1)/(T(V_T^{S^{\star}}+V_T^{\bx^{\star}}))$ and $ R(1)/(T^{2/3}(V_T^{S^{\star}}+V_T^{\bx^{\star}}))$, where the input parameters to Algorithm~\ref{alg:OMDCO H >1} are set as Theorem~\ref{thm:regret of A_3}(b) and $R(1)$ is the $1$-regret of Algorithm~\ref{alg:OMDCO H >1} defined in Eq.~\eqref{eqn:def of R alpha}.  Recall from Theorem~\ref{thm:regret of A_3}(b) that $R(1)\le\tilde{O}((\E[V_T^{S^{\star}}]+\E[V_T^{\bx^{\star}}])T^{2/3})$, which implies $R(1)/(T(\E[V_T^{S^{\star}}]+\E[V_T^{\bx^{\star}}]))\le\tilde{O}(T^{-1/3})$ and $R(1)/(T^{2/3}(\E[V_T^{S^{\star}}]+\E[V_T^{\bx^{\star}}]))\le\tilde{O}(1)$. Since $ R(1)/(T(V_T^{S^{\star}}+V_T^{\bx^{\star}}))$ in Fig.~\ref{fig:dynamic regret}(a) decreases as $T$ increases, Fig.~\ref{fig:dynamic regret}(a) matches with the result in Theorem~\ref{thm:regret of A_3}(b). However, $ R(1)/(T^{2/3}(V_T^{S^{\star}}+V_T^{\bx^{\star}}))$ in Fig.~\ref{fig:dynamic regret}(b) also decreases as $T$ increases, which implies that the regret bound in Theorem~\ref{thm:regret of A_3}(b) may not be tight. In fact, the factor $\E[V_T^{S^{\star}}]+\E[V_T^{\bx^{\star}}]$ can potentially be improved as we argued in Sections~\ref{sec:OMDCO with I=[n]}-\ref{sec:OMDCO with general I}. To obtain the results in Fig.~\ref{fig:dynamic regret limited}, we consider the case when the adversary only switches the function $f_t^s(\cdot,\cdot)$ for at most $\lceil T^{1/6}\rceil$ times, which implies $\E[V_T^{S^{\star}}]+\E[V_T^{\bx^{\star}}]\le O(T^{1/6})$. 
In Fig.~\ref{fig:dynamic regret limited}, we plot $R(1)/(T\ln T)$ and $R(1)/(T^{3/4}\ln T)$, where we set $\gamma=\frac{nT^{1/12}}{T^{1/3}}$ and set the other input parameters to Algorithm~\ref{alg:OMDCO H >1} as Theorem~\ref{thm:regret of A_3}(b). Since  $\E[V_T^{S^{\star}}]+\E[V_T^{\bx^{\star}}]\le O(T^{1/6})$, we obtain from Theorem~\ref{thm:regret of A_3}(b) and our discussions in Section~\ref{sec:OMDCO with general I} that the regret bound becomes $R(1)\le O(T^{3/4}\ln T)$, which implies $R(1)/(T\ln T)\le O(T^{-1/4})$ and $R(1)/(T^{3/4}\ln T)\le O(1)$. Since  $R(1)/(T\ln T)$ in Fig.~\ref{fig:dynamic regret limited}(a) decreases as $T$ increases and $R(1)/(T^{3/4}\ln T)$ remains almost unchanged, the results in Fig.~\ref{fig:dynamic regret limited}(a)-(b) match with the result in Theorem~\ref{thm:regret of A_3}(b) and show that the upper bound $O(T^{3/4}\ln T)$ is tight for the instances of the OSS problem considered above.

\section{Conclusion}
We formulated and studied OMDCO problem, where the decision maker simultaneously chooses a discrete and a continuous variable over $T$ rounds. We proposed algorithms to solve this problem under different settings and proved that the algorithms enjoy sublinear regret in $T$. Our regret analysis extends the existing results on MAB algorithms and OCO algorithms. To demonstrate the applicability of the problem setup, we showed that many applications fit into this setup and the proposed algorithms can be applied to solve these applications with regret guarantees. We validated our results with experiments. Future work includes considering more general reward functions and proving regret lower bounds.

\bibliographystyle{plain}        
\bibliography{autosam}           



\appendix
\section{Proofs in Section~\ref{sec:preliminary results}}\label{app:proof of preliminary results}
\subsection{Proof of Proposition~\ref{prop:MAB with error}}
For any $j\in[n]$ and any $t\in[T]$, we denote 
\begin{equation*}
\tilde{r}_t^j=\frac{\beta_t(r_t^j-\epsilon_t^j+\epsilon_t-a)}{b-a },
\end{equation*}
which yields $\tilde{r}_t^j\in[0,1]$ and $\hat{r}_t^j=\tilde{r}_t^j/p_t^j$ if $j=i_t$.  By the assumption made in Proposition~\ref{prop:MAB with error}, we have $\tilde{r}_t^j\in[0,1]$. Based on this, one can first follow the steps in the proof of \cite[Theorem~3.1]{auer2002nonstochastic} and obtain that 
	\begin{equation*}
		\frac{W_{t+1}}{W_t}\le 1+\frac{\gamma/n}{1-\gamma}\tilde{r}_t^{i_t}+\frac{(e-2)(\gamma/K)^2}{1-\gamma}\sum_{j=1}^n\hat{r}_t^{i_t}+\frac{e}{T},
	\end{equation*}
	where $W_t\triangleq\sum_{j=1}^n\varpi_t^j$. We partition the time horizon $[1,\dots,T]$ into segments $[T_1,\dots,T_2),\dots,[T_{V_T^{i^{\star}}},\dots,T_{V_T^{i^{\star}}+1})$, where $T_1=1$, $T_{V_T^{i^{\star}}+1}=T+1$, and $i^{\star}_{T_s}=\cdots=i^{\star}_{T_{s+1}-1}$ for all $s\in[V_T^{i^{\star}}]$. One can now follow the steps in the proof of \cite[Theorem~8.1]{auer2002nonstochastic} and obtain that 
	\begin{multline*}
		\sum_{t=T_s}^{T_{s+1}-1}\tilde{r}_t^{i_t}\ge(1-\gamma)\sum_{t=T_s}^{T_{s+1}-1}\hat{r}_t^{i_t^{\star}}-\frac{n\ln(nT)}{\gamma}\\-(e-2)\frac{\gamma}{n}\sum_{t=T_s}^{T_{s+1}-1}\sum_{j=1}^n\hat{r}_t^j-\frac{en(T_{s+1}-T_s)}{T\gamma}.
	\end{multline*}
	Summing the above display over all $s\in[V_T^{i^{\star}}]$ yields
\begin{multline}\label{eqn:result from exp3.S}
\sum_{t=1}^T\tilde{r}_t^{i_t}\ge(1-\gamma)\sum_{t=1}^T\hat{r}_t^{i_t^{\star}}-\frac{n(V_T^{i^{\star}} \ln(nT)+e)}{\gamma}\\-(e-2)\frac{\gamma}{n}\sum_{t=1}^T\sum_{i=1}^n\hat{r}_t^i.
\end{multline}
Considering any $j\in[n]$, we have from line~7 of Algorithm~\ref{alg:exp3.S with error} that
\begin{align}\nonumber
\E[\hat{r}^j_t]&=\E\big[\E[\hat{r}^j_t|i_1,\dots,i_{t-1}]\big]\\\nonumber
&=\rho\E\big[p_t^j\frac{r_t^j-\epsilon_t^j+\epsilon_t-a}{p_t^j(b-a)}+(1-p_t^j)0\big]\\\nonumber
&=\frac{\rho}{b-a}\E[r_t^j-\epsilon_t^j+\epsilon_t-a]\\\nonumber
&\le\frac{\rho}{b-a}\big(\E[r_t^j]-\underline{\epsilon}+\E[\epsilon_t]-a\big),
\end{align}
where we use the fact that $\beta_t\sim\textbf{Bernoulli}(\rho)$ is assumed to be independent of the random choice in line~4 of Algorithm~\ref{alg:exp3.S with error}. It then follows that
\begin{multline}\label{eqn:upper bound on r_hat}
\E\Big[\sum_{t=1}^T\sum_{i=1}^n\hat{r}_t^i\Big]\\\le\frac{\rho}{b-a}\big(nbT-\underline{\epsilon}nT+n\E\Big[\sum_{t=1}^T\epsilon_t\Big]-naT\big),
\end{multline}
where we use the fact that $r_t^j\in[a,b]$ for all $t\in[T]$ and all $j\in[n]$. Alternatively, we may write $\E[\hat{r}_t^j]$ as
\begin{align*}
\E[\hat{r}_t^j]&=\rho\E\big[\BI\{i_t=j\}\frac{r_t^j-\epsilon_t^j+\epsilon_t-a}{p_t^j(b-a)}\big]\\
&=\frac{\rho}{b-a}\big(\E[r_t^j+\epsilon_t-a]-\E\big[\BI\{i_t=j\}\frac{\epsilon_{t}^{i_t}}{p_t^{i_t}}\big]\big)\\
&\ge\frac{\rho}{b-a}\big(\E[r_t^j+\epsilon_t]-a-\frac{n}{\gamma}\E[\bar{\epsilon}_{t}^{i_t}]\big),
\end{align*}
where we use the facts that $\BI\{i_t=j\}\le1$ and $1\le 1/p_t^{i_t}\le n/\gamma$. Since $\gamma\in(0,1]$, we further obtain
\begin{multline}\label{eqn:lower bound on r_hat}
\E\Big[(1-\gamma)\sum_{t=1}^T\hat{r}_t^{i_t^{\star}}\Big]\ge \frac{\rho}{b-a}\Big(\E\Big[\sum_{t=1}^Tr_t^{i_t^{\star}}\Big]-\gamma bT\\+(1-\gamma)\E\Big[\sum_{t=1}^T\epsilon_t\Big]-(1-\gamma)\frac{n}{\gamma}\E\Big[\sum_{t=1}^T\bar{\epsilon}_t^{i_t}\Big]-(1-\gamma)aT\Big).
\end{multline}
Finally, taking expectation on both sides of \eqref{eqn:result from exp3.S}, using \eqref{eqn:upper bound on r_hat}-\eqref{eqn:lower bound on r_hat} and rearranging, we obtain \eqref{eqn:regret bound on exp3.S with error}.$\hfill\blacksquare$

\subsection{Proof of Proposition~\ref{prop:OCO with one-point evaluation}}
First, since $\bz_1\in(1-\xi)\X$, we have from \cite[Observation~3.2]{flaxman2005online} and our choices of $\delta$ and $\xi$ that $\bx_t\in\X$ for all $t\in[T]$. Let $\bz_t^{\star}=(1-\xi)\bx_t^{\star}$ for all $t\in[T]$, which yields $\bz_t^{\star}\in(1-\xi)\X$. For our analysis in this proof, we define
\begin{align}
\hat{h}_t(\bz)&=\E_{\bv}[h_t(\bz+\delta\bv)],\label{eqn:def of f_t hat}\\
l_{t}(\bz)&=\hat{h}_{t}(\bz)+(\tilde{\bg}_t-\nabla\hat{h}_{t}\label{eqn:def of l_t}(\bz_{t}))^{\top}\bz,
\end{align}
for all $t\in[T]$ and all $\bz\in(1-\xi)\X$, where $\E_{\bv}[\cdot]$ denotes the expectation with respect to the uniform distribution of $\bv$ over the unit ball $\BB$. Here, we set $\tilde{\bg}_t\triangleq\frac{d}{\delta}h_t(\bx_t)\bu_t$ if the single point feedback is used in Algorithm~\ref{alg:bandit OCO with bandit feedback}, and set $\tilde{\bg}_t\triangleq\frac{d}{2\delta}(h_t(\bx_t)-h_t(\tilde{\bx}_t))\bu_t$ if the two-point feedback is used in Algorithm~\ref{alg:bandit OCO with bandit feedback}. 

\textit{Proof of (a)}: We begin by considering a fictitious round $t=0$ in Algorithm~\ref{alg:bandit OCO with bandit feedback}, where we instead initialize with $\bz_0\in(1-\xi)\X$ and define an auxiliary function $h_0(\bx)=-\frac{T^{3/4}}{2}\norm{\bx}^2$ for $\bx\in\X$. One can show that $h_0(\cdot)$ is $\mu_0$-strongly concave with $\mu_0=T^{3/4}$, and $h_0(\cdot)$ is $T^{3/4}D$-Lipschitz continuous. Similarly, we can define $\hat{h}_0(\cdot)$ and $h_0(\cdot)$, where we set $\tilde{\bg}_0=\frac{d}{2\delta}(h_0(\bx_0)-h_0(\tilde{\bx}_0))\bu_0$ with $\bx_0$ and $\tilde{\bx}_0$ defined similarly as lines~4 in Algorithm~\ref{alg:bandit OCO with bandit feedback}. From the above discussions, we see that $\nabla l_t(\bz)=\nabla\hat{h}_t(\bz)+\tilde{\bg}_t-\nabla\hat{h}_t(\bz_t)$ for all $t\in\{0,\dots,T\}$, which implies that $\nabla l_t(\bz_t)=\tilde{\bg}_t$. Hence, one can view that the sequence $\bz_0,\bz_1,\dots,\bz_T$ is obtained by applying the deterministic online gradient ascent algorithm (e.g., \cite{hazan2007logarithmic}) to the functions $l_0(\cdot),\dots,l_T(\cdot)$ with the update rule $\bz_{t+1}=\Pi_{(1-\xi)\X}(\bz_t+\eta_t\nabla l_t(\bz_t))$. 

To proceed, we first upper bound the dynamic regret of the online gradient ascent algorithm that we just described, which is given by $\sum_{t=0}^Tl_t(\bz_t^{\star})-\sum_{t=0}^Tl_t(\bz_t)$, where $\bz_0^{\star}=(1-\xi)\bx_0^{\star}$ with $\bx_0^{\star}=\bx_1^{\star}$. Recall that $h_0(\cdot)$ is $\mu_0$-strongly concave with $\mu_0=T^{3/4}$ and $h_t(\cdot)$ is concave for all $t\in[T]$ (i.e., $h_t(\cdot)$ is $\mu_t$-strongly concave with $\mu_t=0$ for all $t\in[T]$). From the definitions of $\hat{h}_t(\cdot)$ and $l_t(\cdot)$, one can also show that $\hat{h}_t(\cdot)$ and $l_t(\cdot)$ are $\mu_t$-strongly concave for all $t\in\{0,\dots,T\}$. Since we set $\eta_t=\frac{1}{T^{3/4}}$ for all $t\in\{0,\dots,T\}$ in Algorithm~\ref{alg:bandit OCO with bandit feedback}, it follows that $\eta_t=\frac{1}{\mu_{0:t}}$, where $\mu_{0:t}\triangleq\sum_{k=0}^t\mu_k$. Now, considering any $t\in\{0,\dots,T\}$, we have from the strongly concavity of $l_t(\cdot)$ that 
\begin{align}\nonumber
l_t(\bz_t^{\star})-l_t(\bz_t)\le\nabla l_t(\bz_t)^{\top}(\bz_t^{\star}-\bz_t)-\frac{\mu_t}{2}\norm{\bz_t^{\star}-\bz_t}^2.
\end{align}
Moreover, using the properties of projections onto convex sets (see, e.g., \cite[Lemma~8]{hazan2007logarithmic}), we have
\begin{align}\nonumber
\norm{\bz_t^{\star}-\bz_{t+1}}^2&=\norm{\bz_t^{\star}-\Pi_{(1-\xi)\X}(\bz_{t}+\eta_t\nabla l_t(\bz_t)}^2\\\nonumber
&\le\norm{\bz_t^{\star}-\bz_t-\eta_t\nabla l_t(\bz_t)}^2,
\end{align}
which implies that
\begin{multline*}
\nabla l_t(\bz_t)^{\top}(\bz_t^{\star}-\bz_t)\le\frac{\norm{\bz_t^{\star}-\bz_t}^2-\norm{\bz_t^{\star}-\bz_{t+1}}^2}{2\eta_t}\\+\frac{\eta_t}{2}\norm{\nabla l_t(\bz_t)}^2.
\end{multline*}
Note that 
\begin{align*}
&\norm{\bz_t^{\star}-\bz_{t+1}}^2\\
=&\norm{\bz_t^{\star}-\bz_{t+1}^{\star}+\bz_{t+1}^{\star}-\bz_{t+1}}^2\\
=&\norm{\bz_t^{\star}-\bz_{t+1}^{\star}}^2+\norm{\bz_{t+1}^{\star}-\bz_{t+1}}^2\\
&\qquad\qquad\qquad+2(\bz_t^{\star}-\bz_{t+1}^{\star})^{\top}(\bz_{t+1}^{\star}-\bz_{t+1})\\
\ge&\norm{\bz_{t+1}^{\star}-\bz_{t+1}}^2+\norm{\bz_{t}^{\star}-\bz_{t+1}^{\star}}^2-4D\norm{\bz_{t}^{\star}-\bz_{t+1}^{\star}}\\
\ge&\norm{\bz_{t+1}^{\star}-\bz_{t+1}}^2-4D\norm{\bz_{t}^{\star}-\bz_{t+1}^{\star}}
\end{align*}
where we use the fact that $\bz_{t+1}^{\star},\bz_{t+1}\in(1-\xi)\X\subseteq(1-\xi)D\BB$ which implies that $\norm{\bz_{t+1}^{\star}-\bz_{t+1}}\le2D$. Combining the above arguments together, we have
\begin{align}\nonumber
&l_t(\bz_t^{\star})-l_t(\bz_t)\\\nonumber
\le &\frac{\norm{\bz_t^{\star}-\bz_t}^2-\norm{\bz_{t+1}^{\star}-\bz_{t+1}}^2+4D\norm{\bz_t^{\star}-\bz_{t+1}^{\star}}}{2\eta_t}\\\nonumber
&\qquad\qquad+\frac{\eta_t}{2}\norm{\nabla l_t(\bz_t)}^2-\frac{\mu_t}{2}\norm{\bz_t^{\star}-\bz_{t}}^2.
\end{align}
Summing from $t=0$ to $t=T$ and rearranging terms, we obtain
\begin{align}\nonumber
&\sum_{t=0}^T\big(l_t(\bz_t^{\star})-l_t(\bz_t)\big)\\\nonumber
\le &\sum_{t=1}^T\norm{\bz_t^{\star}-\bz_t}(\frac{1}{2\eta}-\frac{1}{2\eta_{t-1}}-\frac{\mu_t}{2})+\sum_{t=0}^T\frac{\eta_t}{2}\norm{\nabla l_t(\bz_t)}^2\\\nonumber
&\qquad\qquad\qquad\qquad\qquad+\sum_{t=0}^T\frac{2D\norm{\bz_t^{\star}-\bz_{t+1}^{\star}}}{\eta_t}\\
=&\sum_{t=0}^T\frac{\eta_t}{2}\norm{\nabla l_t(\bz_t)}^2+\sum_{t=0}^{T-1}\frac{2D\norm{\bz_t^{\star}-\bz_{t+1}^{\star}}}{\eta_t},\label{eqn:regret bound for h_t}
\end{align} 
where we use the fact that $\eta_t=\frac{1}{\mu_{0:t}}$ for all $t\in\{0,\dots,T\}$ and we let $\bz_{T+1}^{\star}=\bz_T^{\star}$ since $l_{T+1}(\cdot)$ is of no interest.

Now, based on the definition of $l_t(\cdot)$, we can bound $\norm{\nabla l_t(\bz_t)}=\norm{\tilde{\bg}_t}\le\frac{d}{\delta}|h_t(\bx_t)|\norm{\bu_t}\le\frac{dC}{\delta}$ for all $t\in[T]$. Similarly, one can show that $\norm{\nabla l_0(\bz_0)}=\norm{\tilde{\bg}_0}\le\frac{d}{2\delta}|h_0(\bz_0+\delta)-h_0(\bz_0-\delta\bu_0)|\norm{\bu_0}\le T^{3/4}Dd$, where the second inequality follows from the fact that $h_0(\cdot)$ is $T^{3/4}D$-Lipschitz. We then see from \eqref{eqn:regret bound for h_t} that 
\begin{align}\nonumber
&\sum_{t=0}^T\big(l_t(\bz_t^{\star})-l_t(\bz_t)\big)\\\nonumber
\le&\frac{T^{3/2}d^2D^2}{2T^{3/4}}+\frac{d^2C^2}{2\delta^2}\sum_{t=1}^T\frac{1}{T^{3/4}}+\sum_{t=0}^{T-1}\frac{2D\norm{\bz_t^{\star}-\bz_{t+1}^{\star}}}{\eta_t}\\
\le&\frac{d^2(D^2+C^2)r^2}{2}T^{3/4}+\sum_{t=0}^{T-1}\frac{2D\norm{\bz_t^{\star}-\bz_{t+1}^{\star}}}{\eta_t},\label{eqn:regret bound for h_t 2nd}
\end{align}
 where we plug in $\delta=\frac{r}{T^{1/4}}$. Also note from \cite[Lemma~2.1]{flaxman2005online} that for any $t\in\{0,\dots,T\}$, $\E_{\bu_t}[\tilde{\bg}_t]=\nabla\hat{h}(\bz_t)$, which implies that $\E_{\bu_t}[l_t(\bz)]=\hat{h}_t(\bz)$ (for all $\bz$ that is independent of $\bu_t$). Taking the expectation on both sides of \eqref{eqn:regret bound for h_t 2nd} with respect to the randomness of Algorithm~\ref{alg:bandit OCO with bandit feedback} and recalling that $\bz_t^{\star}=(1-\xi)\bx_t^{\star}$ and $\bx_0^{\star}=\bx_1^{\star}$, we get  
\begin{multline}\label{eqn:eqn:result from adaptive OGD f_t hat}
\E\Big[\sum_{t=0}^T\hat{h}_t((1-\xi)\bx_t^{\star})-\sum_{t=0}^T\hat{h}_t(\bz_t)\Big]\le \frac{d^2(D^2+C^2)}{2}T^{3/4}\\+\sum_{t=1}^{T-1}\frac{2D\E\norm{\bx_t^{\star}-\bx_{t+1}^{\star}}}{\eta_t}.
\end{multline}
Again, recalling from our assumption that $h_t(\cdot)$ is $G$-Lipschitz for all $t\in[T]$ and $\X\subseteq D\BB$, and noting that $\bz_t=\bx_t-\delta\bu_t$, one can show that for any $t\in[T]$,
\begin{align}
&\hat{h}_t(\bz_t)\le h_t(\bx_t)+2G\delta,\label{eqn:inter ineqs 1}\\
&\hat{h}_t((1-\xi)\bx_t^{\star})\ge h_t(\bx_t^{\star})-G\delta-G\xi D.\label{eqn:inter ineqs 2}
\end{align}

Finally, let us consider $t=0$. Since $\bz_0^{\star}\in(1-\xi)\X$, where $\X$ is convex with $r\BB\subseteq\X$, and $\xi=\frac{\delta}{r}$, we know from \cite[Observation~3.2]{flaxman2005online} that $(\bz_0^{\star}+\delta\bv)\in\X\subseteq D\BB$. We then have from the definitions of $h_0(\cdot)$ and $\hat{h}_0(\cdot)$ that 
\begin{align}\nonumber
\hat{h}_0(\bz_0^{\star})-\hat{h}_0(\bz_0)\ge-\frac{T^{3/4}}{2}\E_{\bv}\big[\norm{\bz_0^{\star}+\delta\bv}^2\big]\ge-\frac{T^{3/4}D^2}{2}.
\end{align}
Combining the above arguments together and plugging in $\delta=\frac{r}{T^{1/4}}$, $\xi=\frac{\delta}{r}$ and $\eta_t=\frac{1}{T^{3/4}}$, we conclude that
\begin{align}\nonumber
&\E\Big[\sum_{t=1}^Th_t(\bx_t^{\star})-\sum_{t=1}^Th_t(\bx_t)\Big]\\\nonumber
\le& \frac{d^2(D^2+C^2)r^2}{2}T^{3/4}+3TG\delta+TG\xi D\\\nonumber
&\qquad\qquad\qquad\qquad+\frac{T^{3/4}D^2}{2}+2DT^{3/4}\E[V_T^{\bx^{\star}}]\\\nonumber
\le&(\frac{d^2(D^2+C^2)r^2+D^2}{2}+3Gr+\frac{GD}{r}+2D\E[V_T^{\bx^{\star}}])T^{3/4}.
\end{align}

\textit{Proof of (b)}: The proof follows from similar arguments to those in the proof of single point feedback. Again, we consider a fictitious round $t=0$ in Algorithm~\ref{alg:bandit OCO with bandit feedback}, where we instead initialize with $\bz_0\in(1-\xi)\X$ and define an auxiliary function $h_0(\bx)=-\frac{T^{1/2}}{2}\norm{\bx}^2$ for $\bx\in\X$. One can show that $h_0(\cdot)$ is $\mu_0$-strongly concave with $\mu_0=T^{1/2}$, and $h_0(\cdot)$ is $T^{1/2}D$-Lipschitz continuous. Similarly, one can show that $\norm{\nabla l_0(\bz_0)}=\norm{\tilde{\bg}_0}\le\frac{d}{2\delta}|h_0(\bz_0+\delta\bu_0)-h_0(\bz_0-\delta\bu_0)|\norm{\bu_0}\le T^{1/2}Dd$. Recalling that $h_t(\cdot)$ is $G$-Lipschitz for all $t\in[T]$, one can also show that $\norm{\nabla l_t(\bz_t)}=\norm{\tilde{\bg}_t}\le\frac{d}{2\delta}|h_t(\bz_t+\delta\bu_t)-h_t(\bz_t-\delta\bu_t)|\norm{\bu_t}\le Gd$. Now, plugging the above upper bounds on $\norm{\nabla l_t(\bz_t)}$ for all $t\in\{0,\dots,T\}$ into \eqref{eqn:regret bound for h_t}, we obtain
\begin{align}\nonumber
&\sum_{t=0}^T\big(l_t(\bz_t^{\star})-l_t(\bz_t)\big)\\\nonumber
\le&\frac{Td^2D^2}{2T^{1/2}}+\frac{G^2d^2}{2}\sum_{t=1}^T\frac{1}{T^{1/2}}+\sum_{t=0}^{T-1}\frac{2D\norm{\bz_t^{\star}-\bz_{t+1}^{\star}}}{\eta_t}\\
\le&\frac{d^2(D^2+G^2)}{2}T^{1/2}+\sum_{t=1}^{T-1}\frac{2D\norm{\bz_t^{\star}-\bz_{t+1}^{\star}}}{\eta_t},\label{eqn:regret bound for h_t 3rd}
\end{align}
which also implies that 
\begin{multline*}
\E\Big[\sum_{t=0}^T\hat{h}_t((1-\xi)\bx_t^{\star})-\sum_{t=0}^T\hat{h}_t(\bz_t)\Big]\le \frac{d^2(D^2+G^2)}{2}T^{1/2}\\+\sum_{t=1}^{T-1}\frac{2D\E\norm{\bx_t^{\star}-\bx_{t+1}^{\star}}}{\eta_t}.
\end{multline*}
Noting that \eqref{eqn:inter ineqs 1}-\eqref{eqn:inter ineqs 2} still hold and that $\hat{h}_0(\bz_0^{\star})-\hat{h}_0(\bz_0)\ge-\frac{T^{1/2}D^2}{2}$, one can combine the above arguments and show that
\begin{align}\nonumber
&\E\Big[\sum_{t=1}^Th_t(\bx_t^{\star})-\sum_{t=1}^Th_t(\bx_t)\Big]\\\nonumber
\le& \frac{d^2(D^2+G^2)}{2}T^{1/2}+3TG\delta+TG\xi D\\\nonumber
&\qquad\qquad\qquad\qquad+\frac{T^{1/2}D^2}{2}+2DT^{1/2}\E[V_T^{\bx^{\star}}]\\\nonumber
\le&(\frac{d^2(D^2+G^2)+D^2}{2}+3Gr+\frac{GD}{r}+2D\E[V_T^{\bx^{\star}}])T^{1/2},
\end{align}
where we plugged in $\delta=\frac{1}{T^{1/2}}$,  $\xi=\frac{\delta}{r}$ and $\eta_t=\frac{1}{T^{1/2}}$.

\textit{Proof of (c)}: Since we further assume that $h_t(\cdot)$ is $\mu_t$-strongly concave with $\mu_t=\mu\in\R_{>0}$ for all $t\in[T]$, one can also show that $h_t(\cdot)$ is $\mu$-strongly concave. Note that we may still write $\eta_t=\frac{1}{\mu_{1:t}}$ with $\mu_{1:t}=\sum_{k=1}^t\mu_k$. Following similar arguments to those leading up to \eqref{eqn:regret bound for h_t}, one can then show that 
\begin{multline}
\sum_{t=1}^T\big(l_t(\bz_t^{\star})-l_t(\bz_t)\big)
\le \sum_{t=1}^T\frac{\eta_t}{2}\norm{\nabla l_t(\bz_t)}^2\\+\sum_{t=1}^{T-1}\frac{2D\norm{\bz_t^{\star}-\bz_{t+1}^{\star}}}{\eta_t}.\label{eqn:regret bound for h_t 4th}
\end{multline} 
Moreover, we still have the bounds $\norm{\nabla l_t(\bz_t)}\le Gd$ for all $t\in[T]$. Thus, we get 
\begin{align}\nonumber
&\sum_{t=1}^T\big(l_t(\bz_t^{\star})-l_t(\bz_t)\big)\\\nonumber
\le &\frac{\mu G^2d^2}{2}\sum_{t=1}^T\frac{1}{t}+\sum_{t=1}^{T-1}\frac{2D\norm{\bz_t^{\star}-\bz_{t+1}^{\star}}}{\eta_t}\\\nonumber
\le &\frac{\mu G^2d^2}{2}(1+\ln T)+D\mu\sqrt{\big(\sum_{t=1}^{T-1}\norm{\bz_t^{\star}-\bz_{t+1}^{\star}}^2\big)\big(\sum_{t=1}^{T-1}\frac{1}{t^2}\big)}\\\nonumber
\le &\frac{\mu G^2d^2}{2}(1+\ln T)+D\mu\sqrt{\big(\sum_{t=1}^{T-1}\norm{\bz_t^{\star}-\bz_{t+1}^{\star}}\big)^2\big(\sum_{t=1}^{T-1}\frac{1}{t}\big)^2}\\\nonumber
\le & \frac{\mu G^2d^2}{2}(1+\ln T)+D\mu(1+\ln T)\sum_{t=1}^{T-1}\norm{\bz_t^{\star}-\bz_{t+1}^{\star}},
\end{align}
where we use the CS inequality to obtain the second inequality. Noting the choices of $\xi$ and $\delta$, the remaining proof of (c) now follows from those of (a) and (b).$\hfill\blacksquare$

\section{Proofs in Section~\ref{sec:OMDCO}}\label{app:main results proofs}
\subsection{Proof of Theorem~\ref{thm:regret of A_2}}
First, we decompose $R(1)$ as
\begin{align}\nonumber
R(1)&=\underbrace{\E\Big[\sum_{t=1}^Tf_t(i^{\star}_t,\bx^{\star}_t)-\sum_{t=1}^Tf_t(i_t,\bx^{\star}_t)\Big]}_{R_{\A_s}^1}\\
&\qquad+\underbrace{\E\Big[\sum_{t=1}^Tf_t(i_t,\bx^{\star}_t)-\sum_{t=1}^Tf_t(i_t,\bx_t)\Big]}_{R_{\A_s}^2},\label{eqn:regret decomposition of R_A2}
\end{align}
where $i^{\star}_t\in\I$ and $\bx^{\star}_t\in\X$ are defined as \eqref{eqn:optimal point}. Note that $R_{\A_s}^1$ (resp., $R_{\A_s}^2$) can be viewed as the regret incurred by the discrete variable $i_t\in\I$ (resp., the continuous variable $\bx_t\in\X$).

Now, recall that Algorithm~\ref{alg:OMDCO H=1} is applied to the instance of the MAB problem, where the set of possible actions is given by $\I$, and the reward of choosing any action $i\in\I$ in any round $t\in[T]$ is given by $f_t(i_t,\bx_t^{\star})$. Moreover, since Algorithm~\ref{alg:OMDCO H=1} feeds back $f_t(i_t,\bx_t)$ as the reward of choosing $i_t$ in round $t\in[T]$, one can view $f_t(i_t,\bx_t)$ as an erroneous version of $f_t(i_t,\bx_t^{\star})$. Following our notations in Section~\ref{sec:multi-armed bandit}, we let $r_t^j=f_t(j,\bx_t^{\star})$, $\epsilon_t^j=f_t(j,\bx_t^{\star})-f_t(j,\bx_t)$ and $\epsilon_t=0$ for all $t\in[T]$ and all $j\in\I$. We may now write $f_t(i_t,\bx_t)=r_t^{i_t}-\epsilon_t^{i_t}$. From Assumption~\ref{ass:function f_t}, we have $r_t^j\in[0,C]$, $\epsilon_t^j\le\bar{\epsilon}_t^j=f_t(j,\hat{\bx}_t^{\star})-f_t(j,\bx_t)$ and $\epsilon_t^j\ge\underline{\epsilon}=-C$ for all $t\in[T]$ and all $j\in[n]$. One can then apply Proposition~\ref{prop:MAB with error} (with $\rho=1$) and obtain that 
\begin{multline*}
R_{\A_s}^1\le C\frac{n\big(\E[V_T^{i^{\star}}]\ln(nT)+e\big)}{\gamma}-\E\Big[\sum_{t=1}^T\epsilon_t^{i_t}\Big]\\+(e-2)\gamma CT+\frac{n}{\gamma}\E\Big[\sum_{t=1}^T\bar{\epsilon}_t^{i_t}\Big]+C(e-1)\gamma T.
\end{multline*}
Since $\E\big[\sum_{t=1}^T\epsilon_t^{i_t}\big]=R_{\A_s}^2$, we see from Eq.~\eqref{eqn:regret decomposition of R_A2} that
\begin{multline}\label{eqn:R_Ae inter}
R(1)\le C\frac{n\big(\E[V_T^{i^{\star}}]\ln(nT)+e\big)}{\gamma}\\+\frac{n}{\gamma}\E\Big[\sum_{t=1}^T\bar{\epsilon}_t^{i_t}\Big]+C(2e-3)\gamma T.
\end{multline}
Thus, it remains to upper bound $\E\big[\sum_{t=1}^T\bar{\epsilon}_t^{i_t}\big]$.

To proceed, supposing the single-point feedback is used in Algorithm~\ref{alg:OMDCO H=1}, we have $\eta_t=\frac{1}{T^{3/4}}$ for all $t\in[T]$, $\delta=\frac{1}{T^{1/4}}$ and $\xi=\frac{\delta}{r}$. One can then view $\E\big[\sum_{t=1}^T\bar{\epsilon}_t^{i_t}\big]$ as the dynamic regret of Algorithm~\ref{alg:bandit OCO with bandit feedback} (defined in Eq.~\eqref{eqn:dynamic regret of OCO}) when Algorithm~\ref{alg:bandit OCO with bandit feedback} is applied to the sequence of functions $f_1(i_1,\cdot),\dots,f_T(i_T,\cdot)$. From the definitions of Algorithm~\ref{alg:exp3.S with error} and Algorithm~\ref{alg:OMDCO H=1}, we know that for any $t\in[T]$, $i_t$ depends on $\bx_1,\dots,\bx_{t-1}$ but does not depend on $\bx_t$. In other words, one can view that the functions $f_1(i_1,\cdot),\dots,f_T(i_T,\cdot)$ are generated by an adaptive adversary as we described in Section~\ref{sec:OCO with bandit feedback}. It follows that the results of Proposition~\ref{prop:OCO with one-point evaluation} can be applied here. Specifically, we have from Proposition~\ref{prop:OCO with one-point evaluation}(a) that 
\begin{equation}
\E\Big[\sum_{t=1}^T\bar{\epsilon}_t^{i_t}\Big]\le\big(P_1+D\E[\hat{V}_T^{\bx^{\star}}]\big)T^{3/4}.
\end{equation}
Setting $\gamma=\min\{1,\sqrt{\frac{n}{T^{1/4}}}\}$  proves part~(a).\footnote{Note that if $\gamma=1$, i.e., $n\ge T^{1/4}$, the regret bound in part~(a) holds trivially since $f_t(\cdot,\cdot)\in[0,C]$ by Assumption~\ref{ass:function f_t}.} The proof of parts~(b) and~(c) follow similarly by applying Proposition~\ref{prop:OCO with one-point evaluation}(b) and~(c), respectively.$\hfill\blacksquare$

\subsection{Proof of Theorem~\ref{thm:regret of A_3}}
For our analysis in this proof, we introduce the following definitions and notations. First, we define $g_t^{\star}(S)=f_t(S,\bx_t^{\star})$ and $g_t(S)=f_t(S,\bx_t)$ for all $S\subseteq[n]$ and all $t\in[T]$, where $\bx^{\star}_t$ is defined in \eqref{eqn:optimal point} and $\bx_t$ is chosen by Algorithm~\ref{alg:OMDCO H >1} in round $t\in[T]$. Next, we augment the original ground set $[n]=\{1,\dots,n\}$ and define  $\tilde{U}=\{(s_1,\dots,s_T):s_t\in[n],t\in[T]\}$. Hence, any $\tilde{s}\in\tilde{U}$ is a tuple with $T$ elements and let $\tilde{s}_t$ be the $t$th element of $\tilde{s}$. We then define a set function $g^{\star}:2^{\tilde{U}}\to\R_{\ge0}$ such that $g^{\star}(\tilde{S})=\sum_{t=1}^Tg_t^{\star}(\tilde{S}_t)$ for all $\tilde{S}\subseteq\tilde{U}$, where $\tilde{S}_t\triangleq\{\tilde{s}_t:\tilde{s}\in\tilde{S}\}$. Since $g_t^{\star}(\cdot)$ is monotone nondecreasing with $g^{\star}_t(\emptyset)=0$ from Assumption~\ref{ass:function f_t zero}, one can check that $g^{\star}(\cdot)$ is monotone nondecreasing with $g^{\star}(\emptyset)=0$. Similarly to Definitions~\ref{def:submodularity ratio}-\ref{def:curvature}, we can define the submodularity ratio and curvature of $g^{\star}(\cdot)$, denoted as $\gamma_g$ and $c_g$, respectively. One can show via Definitions~\ref{def:submodularity ratio}-\ref{def:curvature} that $\kappa_g=\underline{\kappa}$ and $c_g=\overline{c}$. To proceed, we decompose $R(\alpha)$ as
\begin{align}\nonumber
R(\alpha)&=\underbrace{\E\Big[\alpha\sum_{t=1}^Tf_t(S^{\star}_t,\bx^{\star}_t)-\sum_{t=1}^Tf_t(S_t,\bx^{\star}_t)\Big]}_{R(\alpha)^1}\\
&\quad+\underbrace{\E\Big[\sum_{t=1}^Tf_t(S_t,\bx^{\star}_t)-\sum_{t=1}^Tf_t(S_t,\bx_t)\Big]}_{R(\alpha)^2},\label{eqn:decomposition of R_Au}
\end{align}
where $S_t^{\star}\in\I$ and $\bx_t^{\star}\in\X$ are defined in \eqref{eqn:optimal point}. Let us denote $\tilde{S}^{\star}=\{\tilde{s}^{\star 1},\dots,\tilde{s}^{\star H}\}$ with $\tilde{S}^{\star}_t=S^{\star}_t$ for all $t\in[T]$, $\tilde{S}^{\dagger}=\{\tilde{s}^{\dagger 1},\dots,\tilde{s}^{\dagger H}\}$ with $\tilde{S}^{\dagger}_t=S_t$ for all $t\in[T]$, and $\tilde{S}^{\prime l}=\{\tilde{s}^{\prime 1},\dots,\tilde{s}^{\prime l}\}$ with $\tilde{S}^{\prime l}_t=S^{\prime l}$ for all $t\in[T]$ and all $l\in[H]$, where $\tilde{S}^{\prime 0}=\emptyset$.

First, we upper bound $R(\alpha)^1$ in the above decomposition, which corresponds to the regret incurred by $S_1,\dots,S_T\in\I$. To this end, we recall from the definition of Algorithm~\ref{alg:OMDCO H >1} that the algorithm maintains $H$ independent copies of {\bf Exp3.S}, i.e., $\M_1,\dots,\M_H$, where each $\M_l$ is applied to the instance of MAB, with the set of possible actions given by $[n]$ and the reward of any action $i\in[n]$ in round $t\in[T]$ given by $g_t^{\star}(S_{t}^{\prime l-1}\cup\{i\})-g_t^{\star}(S^{\prime l-1}_t)$. Recalling Eq.~\eqref{eqn:dynamic regret}, for any $l\in[H]$ we define the following dynamic regret of the {\bf Exp3.S} subroutine $\M_l$:
\begin{equation*}
R_l=\max_{\tilde{s}\in\tilde{U},V_T^{\tilde{s}}\le V_T^{S^{\star}}}g^{\star}(\tilde{S}^{\prime l-1}\cup\{\tilde{s}\})-g_t^{\star}(\tilde{S}^{\prime l-1}\cup\{\tilde{s}^l\}),
\end{equation*}
where $V_T^{\tilde{s}}\triangleq 1+\sum_{t=1}^{T-1}\BI\{\tilde{s}_t\neq\tilde{s}_{t+1}\}$, and $\tilde{s}^l=(s_1^l.\dots,s_T^l)$. Note that the maximization is over $\tilde{s}\in\tilde{U}$ such that $V_T^{\tilde{s}}\le V_T^{S^{\star}}$, and one can check that $\tilde{S}^{\star}$ described above satisfies that $V^{\tilde{s}^{\star l}}_T\le V_T^{S^{\star}}$ for all $l\in[H]$, where $V_T^{\tilde{s}^l}$ is defined in the same way as $V_T^{\tilde{s}}$ defined above. It then follows from Lemma~\ref{lemma:greedy algorithm} in Appendix~\ref{app:aux results} that 
\begin{align}\nonumber
g^{\star}(S^{\prime H})&\ge\frac{1}{\overline{c}}(1-e^{-\overline{c}\underline{\kappa}})\max_{\substack{\tilde{S}\subseteq\tilde{U},|\tilde{S}|\le H\\V_T^{\tilde{s}^l}\le V_T^{S^{\star}},\forall l\in[H]}}g^{\star}(\tilde{S})-\sum_{l=1}^H R_l\\
&\ge\frac{1}{\overline{c}}(1-e^{-\overline{c}\underline{\kappa}})g^{\star}(\tilde{S}^{\star})-\sum_{l=1}^H R_l,\label{eqn:greedy solution}
\end{align}
where $\tilde{S}=\{\tilde{s}^1,\dots,\tilde{s}^H\}$ and $V_T^{\tilde{s}^l}$ is defined in the same way as $V_T^{\tilde{s}}$ defined above. Going back to $R_{\A_u}^1$, we get
\begin{align}\nonumber
R_{\A_u}^1&=\E\Big[\alpha g^{\star}(\tilde{S}^{\star})-g^{\star}(\tilde{S}^{\prime H})+g^{\star}(\tilde{S}^{\prime H})-g^{\star}(\tilde{S}^{\dagger})\Big]\\\nonumber
&\le\E\Big[\sum_{l=1}^H R_l\Big]+\E\Big[\sum_{t=1}^T\big(g_t^{\star}(S^{\prime H}_t)-g^{\star}_t(S_t^{\dagger})\big)\Big]\\\nonumber
&=\E\Big[\sum_{l=1}^H R_l\Big]+\E\Big[\sum_{t=1}^T\BI\{\tilde{\beta}_t=1\}C\Big]\\
&=\E\Big[\sum_{l=1}^H R_l\Big]+\tilde{\rho} C T,\label{eqn:R_Au^1 upper bound 1}
\end{align}
where we again use the definition of Algorithm~\ref{alg:OMDCO H >1}.

We then aim to upper bound $\sum_{l=1}^H\E[R_l]$. Note that for any $t\in[T]$ and any $l\in[H]$, Algorithm~\ref{alg:OMDCO H >1} feeds back $\BI\{\tilde{\beta}_t=1,l_t=l,i_t=s_t^l\}g_t(S_t)$ to $\M_l$ as the reward of choosing $s_t^l\in[n]$ in round $t\in[T]$, where we observe from line~7 in Algorithm~\ref{alg:OMDCO H >1} that $\BI\{\tilde{\beta}_t=1,l_t=l,i_t=s_t^l\}g_t(S_t)=\BI\{\tilde{\beta}_t=1,l_t=l,i_t=s_t^l\}g_t(S_t^{\prime l_t-1}\cup\{s_t^l\})$. Following the notations in Section~\ref{sec:multi-armed bandit} we let $r_{t,l}^j=g_t^{\star}(S_t^{\prime l-1}\cup\{j\})-g_t^{\star}(S_t^{\prime l-1})$, $\epsilon_{t,l}^j=g_t^{\star}(S_t^{\prime l-1}\cup\{j\})-g_t(S_t^{\prime l-1}\cup\{j\})$, and $\epsilon_{t,l}=g_t^{\star}(S_t^{\prime l-1})$ for all $j\in[n]$, all $t\in[T]$ and all $l\in[H]$. We can then write $g_t(S_t^{\prime l-1}\cup\{s_t^l\})=r_{t,l}^{s_t^l}-\epsilon_{t,l}^{s_t^l}+\epsilon_{t,l}$. Note from Assumption~\ref{ass:x_t^star optimal} that $\epsilon_{t,l}^j\ge0$ and note from Assumptions~\ref{ass:function f_t}-\ref{ass:function f_t zero} that $r_{t,l}^j\in[0,C]$, $\epsilon_{t,l}^j\ge\underline{\epsilon}=-C$ and $\epsilon_{t,l}\le C$, for all $j\in[n]$, all $t\in[T]$ and all $l\in[H]$. Also note that $\BI\{\tilde{\beta}_t=1,l_t=l,i_t=s_t^l\}$ for all $t\in[T]$ can be viewed as i.i.d. Bernoulli random variables with parameter $\rho=\frac{\tilde{\rho}}{Hn}$. Based on the above arguments, one can now apply Proposition~\ref{prop:MAB with error} and obtain that for any $l\in[H]$,
\begin{multline}
\E[R_l]\le C\frac{n\big(\E[V_T^{S^{\star}}]\ln(nT)+e\big)}{\gamma\rho}+(e-2)\gamma CT\\+\frac{n}{\gamma}\E\Big[\sum_{t=1}^T\epsilon_{t,l}^{s^l_t}\Big]+2(e-1)\gamma CT.
\end{multline}
Hence, to upper bound $\sum_{l=1}^H\E[R_l]$, we need to upper bound $\E\big[\sum_{t=1}^T\epsilon_{t,l}^{s_t^l}\big]$. We have the following:
\begin{align}\nonumber
&\sum_{l=1}^H\E\Big[\sum_{t=1}^T\epsilon_{t,l}^{s_t^l}\Big]\\\nonumber
=&\frac{Hn}{\tilde{\rho}}\sum_{l=1}^H\E\Big[\sum_{t=1}^T\BI\{\tilde{\beta}_t=1,l_t=l,i_t=s_t^l\}\\\nonumber
&\qquad\times\big(g^{\star}_t(S^{\prime l-1}_t\cup\{s_t^l\})-g_t(S^{\prime l-1_t}\cup\{s_t^l\})\big)\Big]\\\nonumber
=&\frac{Hn}{\rho}\E\Big[\sum_{t=1}^T\BI\{\tilde{\beta}_t=1,i_t=s_t^l\}\\\nonumber
&\qquad\times\big(g^{\star}_t(S^{\prime l-1}_t\cup\{s_t^l\})-g_t(S^{\prime l-1_t}\cup\{s_t^l\})\big)\Big]\\
=&\frac{Hn}{\tilde{\rho}}\E\Big[\sum_{t=1}^T\big(g^{\star}_t(S_t)-g_t(S_t)\big)\Big]=\frac{Hn}{\tilde{\rho}}R_{\A_u}^2,\label{eqn:epsilon in alg4}
\end{align}
where the first equation follows from the fact that $\tilde{\beta}_t$, $l_t$ and $i_t$ are independent, and the third equation follows from the definition of Algorithm~\ref{alg:OMDCO H >1}.

Therefore, combining \eqref{eqn:decomposition of R_Au} and \eqref{eqn:R_Au^1 upper bound 1}-\eqref{eqn:epsilon in alg4}, we obtain 
\begin{multline}
R(\alpha)\le C\frac{n\big(\E[V_T^{S^{\star}}]\ln(nT)+e\big)}{\gamma\tilde{\rho}}+(3e-4)\gamma CT\\+2\frac{Hn^2}{\gamma\tilde{\rho}}R_{\A_u}^2+\tilde{\rho} CT.
\end{multline}
Now, following similar arguments to those in the proof of Theorem~\ref{thm:regret of A_2}, one can apply Proposition~\ref{prop:OCO with one-point evaluation}(b)-(c) and obtain upper bounds on $R_{\A_u}^2$. Finally, plugging the choices of $\gamma$ and $\tilde{\rho}$ and using similar  arguments to those in Footnote~11 complete the proof of the theorem.$\hfill\blacksquare$

\section{Proof in Section~\ref{sec:applications}}\label{app:application proof}
\subsection{Proof of Proposition~\ref{prop:sub ratio and curvature in SS 2}}
First, consider any $S,\Omega\subseteq[n]$. We have from the concavity of $f_t^s(S_0,\cdot)$ for any $S_0\subseteq[n]$ that 
\begin{align}\nonumber
&f_t^s(\bx_t^{\star}(S\cup\Omega))-f_t^s(\bx_t^{\star}(S))\\\nonumber
\le&\nabla f_t^{s}(\bx_t^{\star}(S))^{\top}\bx_t^{\star}(S\cup\Omega\setminus S)\\\nonumber
\le&\sum_{\omega\in\Omega\setminus S}|(\nabla f_t^{s}(\bx_t^{\star}(S)))_{\omega}||(\bx_t^{\star})_{\omega}|.
\end{align}
For any $\omega\in\Omega\setminus S$, define $\tilde{f}_t^s(x)=f_t^s(\bx_t^{\star}(S)+\mathbf{1}(\omega)x)$ for all $x\in\R$, where we drop the dependency of $\tilde{f}_t^s(\cdot)$ on $S$ and $\omega$ for notational simplicity. We then have 
\begin{align}\nonumber
\dot{\tilde{f}}_t^{s}(x)&=\nabla f_t^s(\bx_t^{\star}(S)+\mathbf{1}(\omega)x)^{\top}\mathbf{1}(\omega)\\\nonumber
&=(\nabla f_t(\bx_t^{\star}(S)+\mathbf{1}(\omega)x))_{\omega}.
\end{align}
It follows from the mean value theorem that 
\begin{align}\nonumber
f_t^s(\bx_t^{\star}(S\cup\{\omega\}))-f_t^s(\bx_t^{\star}(S))&=\tilde{f}_t^s((\bx_t^{\star})_{\omega})-\tilde{f}^s_t(\mathbf{0})\\\nonumber
&=((\bx_t^{\star})_{\omega}-0)\dot{\tilde{f}}_t^{s}(x^{\prime}),
\end{align}
for some $x^{\prime}\in\X^{\star}_{\omega}$, where $\X^{\star}_{\omega}=(0,(\bx_t^{\star})_{\omega})$ if $(\bx_t^{\star})_{\omega}>0$ and $\X^{\star}_{\omega}=((\bx_t^{\star})_{\omega},0)$ if $(\bx_t^{\star})_{\omega}<0$. Since $f_t^s(\bx_t^{\star}(S\cup\{\omega\}))-f_t^s(\bx_t^{\star}(S))\ge0$ from Assumption~\ref{ass:function f_t zero}, we get
\begin{multline*}
f_t^s(\bx_t^{\star}(S\cup\{\omega\}))-f_t^s(\bx_t^{\star}(S))\\=|(\nabla f_t(\bx_t^{\star}(S)+\mathbf{1}(\omega)x^{\prime}))_{\omega}||(\bx_t^{\star})_{\omega}|.
\end{multline*}
Combining the above arguments and recalling Definition~\ref{def:submodularity ratio}, we get that 
\begin{align}\nonumber
\kappa_t^s&\ge\min_{S,\Omega\subseteq[n]}\frac{\sum_{\omega\in\Omega\setminus S}|(\nabla f_t(\bx_t^{\star}(S)+\mathbf{1}(\omega)x^{\prime}))_{\omega}||(\bx_t^{\star})_{\omega}|}{\sum_{\omega\in\Omega\setminus S}|(\nabla f_t^{s}(\bx_t^{\star}(S)))_{\omega}||(\bx_t^{\star})_{\omega}|}\\\nonumber
&\ge\min_{S,\Omega\subseteq[n],\omega\in\Omega\setminus S}\frac{|(\nabla f_t(\bx_t^{\star}(S)+\mathbf{1}(\omega)x^{\prime}))_{\omega}|}{|(\nabla f_t^{s}(\bx_t^{\star}(S)))_{\omega}|}\\\nonumber
&\ge\min_{\omega\in[n]}\frac{\min_{\bx\in\X}|(\nabla f_t^s(\bx))_{\omega}|}{\max_{\bx\in\X}|(\nabla f_t^s(\bx))_{\omega}|},
\end{align}
which completes the proof of the proposition.$\hfill\blacksquare$

\section{Auxiliary Result}\label{app:aux results}
\begin{lem}\label{lemma:greedy algorithm}
Consider a monotone nondecreasing set function $g:2^{[n]}\to\R_{\ge0}$ with $g(\emptyset)=0$. Let $\kappa_g,c_g\in[0,1]$ be the submodularity ratio and curvature of $g(\cdot)$ given by Definitions~\ref{def:submodularity ratio}-\ref{def:curvature}, respectively. Let $S^{\dagger}=\{s^{\dagger 1},\dots,s^{\dagger H}\}\subseteq[n]$ and denote $S^{\dagger l}=\{s^{\dagger 1},\dots,s^{\dagger l}\}$ for all $l\in[H]$ with $S^{\dagger 0}=\emptyset$. Suppose that $g(S^{\dagger l-1}\cup\{s^{\dagger l}\})\ge\max_{s\in[n]}g(S^{\dagger l-1}\cup\{s\})-\tau_l$ for all $l\in[H]$, where $\tau_l\in\R_{\ge0}$. Then, 
\begin{equation}
g(S^{\dagger})\ge\frac{1}{c_g}(1-e^{-c_g\kappa_g})g(S^{\star})-\sum_{l=1}^H \tau_l,
\end{equation}
where $S^{\star}\in\argmax_{S\subseteq[n],|S|\le H}g(S)$.
\end{lem}
\begin{pf}
First, let us denote $\varphi_l=g(S^{\dagger l-1}\cup\{s^{\dagger l}\})-g(S^{\dagger l-1})$ for all $l\in[H]$. Following similar arguments to those for \cite[Lemma~1]{bian2017guarantees}, one can show that for any $\Omega\subseteq[n]$ with $|\Omega|=H$ and any $l\in\{0,\dots,H-1\}$,
\begin{multline}
c_h\sum_{m:s^{\dagger m}\in S^{\dagger l}\setminus\Omega}\varphi_m + \sum_{m:s^{\dagger m}\in S^{\dagger l}\setminus\Omega}\varphi_l\\+\frac{H-\nu^l}{\kappa_g}(\varphi_{l+1}+\tau_{l+1})\ge g(\Omega),
\end{multline} 
where $\nu^l=|S^{\dagger l}\cap\Omega|$. Moreover, using similar arguments to those for \cite[Theorem~1]{bian2017guarantees}, one can show that the value of $g(S^{\dagger})/g(S^{\star})$ is lower bounded by the optimal solution to the following linear program:
\begin{equation}\label{eqn:LP}
\begin{split}
&\min\sum_{l=1}^H q_l\\
s.t.\ &q_l\ge0,\forall l\in[H],\\
&\mathbf{L}\mathbf{q}\ge\mathbf{1}+\frac{H}{\kappa_g g(S^{\star})}\bm{\tau},
\end{split}
\end{equation}
where $\mathbf{q}=[q_1\ \cdots\ q_H]^{\top}$, $\bm{\tau}=[\tau_1\ \cdots\ \tau_H]$, and $\mathbf{L}\in\R^{H\times H}$ is the matrix with all the diagonal elements to be $\frac{H}{\kappa_g}$, all the upper diagonal elements to be zero and all the off-diagonal elements to be $c_g
$. Now, one can check that the optimal solution to \eqref{eqn:LP} is given by $q_l=\frac{\kappa_g}{H}\big(\frac{H-\kappa_g c_g}{H}\big)^{l-1}+\frac{\tau_l}{g(S^{\star})}$ for all $l\in[H]$. Hence, by combining the above arguments, one can show
\begin{align}\nonumber
\frac{g(S^{\dagger})}{g(S^{\star})}&\ge\frac{1}{c_g}\Big(1-\big(\frac{H-c_g\kappa_g}{H}\big)^H\Big)-\frac{\sum_{l=1}^H \tau_l}{g(S^{\star})}\\\nonumber
&\ge\frac{1}{c_g}(1-e^{-c_g\kappa_g})-\frac{\sum_{l=1}^H \tau_l}{g(S^{\star})},
\end{align}
which completes the proof of the lemma.$\hfill\blacksquare$.
\end{pf}

\end{document}